\documentclass[a4paper]{article}

\usepackage{amsmath}

\usepackage{amsthm}
\usepackage[dvips]{graphicx}

\usepackage{epsfig}
\usepackage{float}

\usepackage{amssymb}
\usepackage{latexsym}

\makeatletter
\newcommand\footnoteref[1]{\protected@xdef\@thefnmark{\ref{#1}}\@footnotemark}
\makeatother

\newtheorem{Theorem}{Theorem}
\newtheorem{Proposition}{Proposition}
\newtheorem{Lemma}{Lemma}
\newtheorem{Corollary}{Corollary}

\newtheorem{Claim}{Claim}

\newtheorem*{rep@theorem}{\rep@title}
\newcommand{\newreptheorem}[2]{%
\newenvironment{rep#1}[1]{%
 \def\rep@title{#2 \ref{##1}}%
 \begin{rep@theorem}}%
 {\end{rep@theorem}}}
\makeatother

\newreptheorem{theorem}{Theorem}
\newreptheorem{corollary}{Corollary}

\theoremstyle{definition}
\newtheorem{definition}{Definition}
\newtheorem{remark}{Remark}

\def\demo{ {\bf Proof.} }

\def\QQ{{\mathbb Q}}

\def\ZZ{{\mathbb Z}}
\def\CC{{\mathbb C}}

\def\N{{\mathcal N}}
\def\C{{\mathcal C}}
\def\Z{{\mathcal Z}}

\def\O{{\mathcal O}}

\def\S{{\mathbf S}}
\def\D{{\mathbf D}}

\begin{document}

\title{Finite group actions on 3-manifolds and cyclic branched covers of knots}

\author{M. Boileau\footnote{\protect\label{f1}Partially supported by ANR
project 12-BS01-0003-01},
C. Franchi,
M. Mecchia\footnote{\protect\label{f2}Partially supported by the FRA 2015 grant
``Geometria e topologia delle variet\`{a} ed applicazioni'', Universit\`{a} di 
Trieste},
L. Paoluzzi\footnoteref{f1},
and B. Zimmermann\footnoteref{f2}}
\date{\today}
\maketitle

\begin{abstract}
\vskip 2mm

As a consequence of a general result about finite group actions on 
$3$-manifolds, we show that a hyperbolic $3$-manifold can be the cyclic 
branched cover of at most fifteen inequivalent knots in $\S^3$  (in fact, a 
main motivation of the present paper is to establish the existence of such a 
universal bound). A similar, though weaker, result holds for arbitrary 
irreducible $3$-manifolds: an irreducible $3$-manifold can be a cyclic branched 
cover of odd prime order of at most six knots in $\S^3$. We note that in most 
other cases such a universal bound does not exist.

\vskip 2mm


\vskip 2mm


\end{abstract}

\section{Introduction}

A classical way to construct closed, connected, orientable $3$-manifolds is
to consider cyclic covers of the $3$-sphere branched along knots. A natural
question in this setting is to understand in how many different ways a closed,
connected, orientable $3$-manifold can be presented as the (total space of a)
cyclic branched cover of a knot. There is an extensive literature on this 
problem, mainly focussing on the case where the branching order is fixed. For
instance, it is known that a closed hyperbolic $3$-manifold is the $n$-fold
cyclic branched cover of at most two knots in $\S^3$, provided $n>2$ \cite{Z2}, 
and at most nine if $n=2$ \cite{Re}. For arbitrary closed, connected, 
orientable, irreducible $3$-manifolds some results are known but only for prime 
orders of ramification. More precisely, such manifolds can be $2$-fold cyclic 
branched covers of arbitrary many knots \cite{Mon1, Mon2} but can cover at most 
two knots if the order is an odd prime \cite{BoPa}. For branching order equal 
to $2$, several analyses of how the different quotient knots are related can 
also be found (see, for instance, \cite{MR} for hyperbolic manifolds, 
\cite{Mon1, Mon2, V, MonW} for toroidal ones, \cite{P} for Conway reducible 
hyperbolic knots, and \cite{Gr} for alternating ones).  

Possibly due to the fact that the knots that are covered by the same manifold
but with different orders of branching are harder to relate, not much was known
so far on the general problem, even for hyperbolic manifolds (see \cite{RZ}). 
The prime motivation for this work was to understand whether it is possible to 
establish bounds on the number of ways a $3$-manifold can be presented as a 
cyclic branched cover of a knot, without fixing the order of the cover. By the
above discussion, no universal bound can be given in general, however the first
main result of the present paper assures the existence of a universal bound for
the class of closed hyperbolic $3$-manifolds.

%A classical, much considered class of closed orientable $3$-manifolds consists 
%of the cyclic branched covers of knots in the $3$-sphere. 
%In the present paper we try to understand in how many different ways a given 
%$3$-manifold can occur as a cyclic branched cover of a knot. For a fixed 
%branching order, the problem  is well-understood (\cite{V}, \cite{MonW}, 
%\cite{Z2}, \cite{Re} and \cite{Gr}); in particular,  a hyperbolic 3-manifold 
%is a 2-fold branched cover of at most nine knots in $\S^3$ and, for $n > 2$, 
%an $n$-fold cyclic branched cover of at most two knots. For different 
%branching orders, the situation turned out to be much more difficult, and even %in the hyperbolic case the existence of a universal bound on the number of 
%different branching orders was not known; the existence of such a universal 
%bound for hyperbolic 3-manifolds is the first main result of the present paper. 

\begin{Theorem}\label{cor:hyperbolic}
A closed hyperbolic $3$-manifold is a cyclic branched cover of at most fifteen 
inequivalent knots in $\S^3$.
\end{Theorem}

We call two knots \emph{equivalent} if one is mapped to the other by an 
orientation-preserving diffeomorphism of $\S^3$. In the present paper, all 
manifolds are closed, connected and orientable, and all maps are smooth and 
orientation-preserving. 

\medskip

The orientation-preserving isometry group of a closed hyperbolic $3$-manifold 
$M$ is finite, and every finite group occurs for some hyperbolic $M$. Suppose 
that $M$ is a cyclic branched cover of a knot in $\S^3$; then the group of 
covering transformations acting on $M$ is generated by a hyperelliptic 
rotation: we call a periodic diffeomorphism of a closed $3$-manifold a 
\emph{hyperelliptic rotation} if all of its non-trivial powers have connected, 
non-empty fixed-point set (a simple closed curve), and its quotient (orbit) 
space is $\S^3$. By the geometrization of $3$-orbifolds, or of finite group 
actions on $3$-manifolds (\cite{BLP}, \cite{BoP}, \cite{DL}), the group of 
covering transformations is conjugate to a subgroup of the isometry group of 
$M$. Hence establishing a universal upper bound for hyperbolic $3$-manifolds as 
in Theorem~\ref{cor:hyperbolic} is equivalent to bounding the number of 
conjugacy classes of cyclic subgroups generated by hyperelliptic rotations of 
the isometry group of a hyperbolic $3$-manifold. Now 
Theorem~\ref{cor:hyperbolic} is a consequence of the following more general 
result on finite group actions on closed $3$-manifolds.

\begin{Theorem}\label{th:main4}
Let $M$ be a closed $3$-manifold not homeomorphic to $\S^3$. Let $G$ be a 
finite group of orientation-preserving diffeomorphisms of $M$. Then $G$ 
contains at most fifteen conjugacy classes of cyclic subgroups generated by a 
hyperelliptic rotation (at most six for cyclic subgroups whose order is not a 
power of two).
\end{Theorem}

Note that the $3$-sphere is the $n$-fold cyclic branched cover of the trivial 
knot for any integer $n\ge2$ (and, by the solution to the Smith conjecture,  
only of the trivial knot). It is well-known that, for any branching order $n$, 
a $3$-manifold can be the $n$-fold cyclic branched cover of an arbitrary number 
of non-prime knots (such a manifold is not irreducible), and that an 
irreducible $3$-manifold can be the $2$-fold branched cover of arbitrarily many 
prime knots (see \cite{Mon1, Mon2, V}). Moreover, Proposition~\ref{prop:circle
bundles} in Section~\ref{sec:seifert} shows that Seifert fibred manifolds can
cover an arbitrary number of knots, all with non-prime orders. It is thus
natural to restrict our attention to covers of odd prime order.

\medskip

For irreducible $3$-manifolds, the following holds:

\begin{Theorem}\label{th:knots}
Let $M$ be a closed, irreducible $3$-manifold. Then there are at most six 
inequivalent knots in $\S^3$ having $M$ as a cyclic branched cover of odd prime 
order.
\end{Theorem}

\medskip

The proof of Theorem~\ref{th:knots} uses Theorem~\ref{th:main4} in connection 
with the equivariant torus-decomposition of irreducible $3$-manifolds into 
geometric pieces, see \cite{BoP}, \cite{BLP}, \cite{CHK}, and \cite{KL}. For 
arbitrary $3$-manifolds, as a direct consequence of Theorem~\ref{th:knots} as 
well as the equivariant prime decomposition for $3$-manifolds \cite{MSY}, the 
following remains true.

\begin{Corollary}\label{th:main}
Let $M$ be a closed $3$-manifold not homeomorphic to $\S^3$. Then $M$ is a 
$p$-fold cyclic branched cover of a knot in $\S^3$ for at most six distinct 
odd prime numbers $p$.
\end{Corollary}

%We note that, in both Theorem~\ref{th:knots} and Corollary~\ref{th:main},  
%there is no universal bound  for  nonprime orders  
%(see Proposition \ref{prop:circle bundles}  in  Section~\ref{sec:seifert}).

\medskip

Another consequence of Theorem~\ref{th:main4} is the following characterisation of the $3$-sphere which generalises the main result of \cite{BPZ} from the case 
of $\ZZ$-homology $3$-spheres to arbitrary closed $3$-manifolds.

\begin{Corollary}\label{cor:char1}
A closed $3$-manifold $M$ is homeomorphic to $\S^3$ if and only if there is a 
finite group $G$ of orientation-preserving diffeomorphisms of $M$ such that $G$ 
contains sixteen conjugacy classes of subgroups generated by hyperelliptic 
rotations. 
\end{Corollary}

A noteworthy aspect of the proof of Theorem~\ref{th:main4} is the substantial 
use of finite group theory, in particular of the classification of finite 
simple groups. For a prime $p$, the $p$-fold cyclic branched cover of a knot in 
$\S^3$ is a rational homology $3$-sphere, and we will prove in 
Section~\ref{sec:nonfree} that every finite group acts non-freely on some 
rational homology $3$-sphere, so the use of the classification seems to be 
intrinsic to the proofs of our results. We note that the class of finite groups 
acting on a $\ZZ/2$-homology $3$-sphere instead is quite restricted (see 
\cite{MZ}), and in this case the much shorter Gorenstein-Harada classification 
of finite simple groups of sectional $2$-rank at most four is sufficient for 
our proofs. The bounds that can be derived for $\ZZ/2$-homology $3$-spheres 
are, however, precisely the same as those we get for arbitrary manifolds.  

\medskip

We have tried to separate the algebraic, purely group theoretical parts of the 
proof (Section~\ref{sec: algebraic lemmas}) from the topological parts 
(Sections~\ref{sec:rotations} and \ref{sec:proof}), so they can be read 
independently.

\medskip

We note that, although the upper bounds in our results are quite small, at this
point we do not know if they are really optimal. In the hyperbolic case, one
can easily construct manifolds that are covers of orders $>2$ of three 
distinct knots \cite{RZ}. For hyperbolic double covers, the bound (nine) is 
sharp according to a result of Kawauchi \cite{Ka}, but no explicit examples are 
known so far. For general irreducible manifolds, Brieskorn spheres of type 
$\Sigma(p,q,r)$, where $p$, $q$, and $r$ are three pairwise different odd 
primes, provide examples of manifolds that cover four knots: $\Sigma(p,q,r)$ 
is the $p$-fold (resp. $q$-fold and $r$-fold) cyclic cover of $\S^3$ branched 
along the torus knot $T(q,r)$ (resp. $T(p,r)$ and $T(p,q)$) as well as the 
double branched cover of a Montesinos knot.

%in fact, a main point of our paper was to establish the existence of such 
%universal upper bounds.  

\medskip

The paper is organised as follows. In Section~\ref{sec:road-map} we present a 
brief sketch of the proof of Theorem~\ref{th:main4}. Hyperelliptic rotations 
and their properties are considered in Section~\ref{sec:rotations}. 
Section~\ref{sec: algebraic lemmas} contains the main group-theoretical part of 
the paper, Section~\ref{sec:proof} the proof of Theorem~\ref{th:main4}, and 
Section~\ref{sec:knots} the proof of Theorem~\ref{th:knots} for the irreducible 
case. Finally, in an Appendix we prove that every finite group acts non-freely 
on some rational homology sphere (adapting the result of \cite{CL} that deals
with free actions).

%%%%%%%%%%%%%%%%%%%%%%%%%%%%%%%%%%%%%%%%%%%%%%%%%%%%%%%%%%%%%%%%%%%%%%%%%% 

\section{Sketch of the proof of Theorem~\ref{th:main4} }\label{sec:road-map}

The proof of Theorem~\ref{th:main4} is based on a series of preliminary
%, sometimes quite technical 
results which are presented in Sections~\ref{sec:rotations} and 
\ref{sec: algebraic lemmas}. Our choice to present the group-theoretical part 
of the proof in a separate section (Section~\ref{sec: algebraic lemmas}) allows for the group-theoretical results to be read independently of the other parts 
of the paper. In the following, in order to make the paper more accessible, 
%readable, 
we explain the main steps of the proof. 

We begin with a more detailed definition of hyperelliptic rotation. 
Note that throughout the paper, unless otherwise stated, \emph{$3$-manifold} 
will mean orientable, connected, closed $3$-manifold. Also, all finite group 
actions by diffeomorphisms will be faithful and orientation-preserving. 

\begin{definition}
Let $\psi:M\longrightarrow M$ be a finite order diffeomorphism of a 
$3$-manifold $M$. We shall say that $\psi$ is a \emph{rotation} if it preserves 
the orientation of $M$, $Fix(\psi)$ is non-empty and connected, and 
$Fix(\psi)=Fix(\psi^k)$ for all non-trivial powers $\psi^k$ of $\psi$. 
$Fix(\psi)$ will be referred to as the \emph{axis of the rotation}. Note that
if $\psi$ is a periodic diffeomorphism of prime order, then $\psi$ is a
rotation if and only if $Fix(\psi)=\S^1$. We shall say that a rotation $\psi$ 
is \emph{hyperelliptic} if the space of orbits $M/\psi$ of its action is 
$\S^3$, and a \emph{hyperelliptic group} is a cyclic group generated by a 
hyperelliptic rotation.
\end{definition}

%We note that, in the situation of Theorem 2, there is a strong dichotomy 
%between the cases of solvable and non-solvable groups $G$. Whereas the case of 
%solvable groups is easy and quite well-understood, it is the non-solvable 
%case that results in the major difficulties since in this case it is much more 
%difficult to relate the Sylow subgroups for different primes. For example, for 
%a fixed, prime branching order, the situation is completely determined by a 
%single Sylow subgroup of $G$, and for more general solvable groups, one can 
%use Hall subgroups instead of Sylow subgroups. 

%In particular, 

We start by observing that the case of hyperelliptic rotations whose order is a 
power of two is already well-understood by work of Reni and Mecchia (see 
\cite{Re} and \cite{Mec0}). In particular, there are at most nine conjugacy
classes of cyclic groups generated by such hyperelliptic rotations. As a
consequence, from now on, we exclude this case and consider only hyperelliptic
groups whose order is not a power of two.

%, and it is not difficult to prove that there are at most nine conjugacy 
%classes of cyclic groups generated by such hyperelliptic rotations (see 
%\cite{Re} and \cite{Mec0}). So in the remaining part of this section, we 
%exclude this case and consider only hyperelliptic groups whose order is not a 
%power of two. 

Section~\ref{sec:rotations} collects various simple facts on the geometry of 
hyperelliptic rotations. In particular, we prove that there are at most 
three hyperelliptic groups commuting pairwise. This result implies the 
existence of a universal bound in the solvable case. In fact, in the solvable 
case, the presence of Hall subgroups assures that all hyperelliptic subgroups 
commute, up to conjugacy, implying that there at most three conjugacy classes 
of such groups.

%As noted above, the most intricate 
The case of non-solvable groups, where local approaches on the basis of 
$p$-groups fail, is more involved. We need a global description of the groups 
that may arise, which is provided in Section~\ref{sec: algebraic lemmas}. In 
that section, we introduce the notion of an \emph{algebraically hyperelliptic 
collection} of cyclic subgroups which have the same algebraic properties as the 
hyperelliptic subgroups; this allows a purely algebraic approach in 
Section~\ref{sec: algebraic lemmas}. The main result there
(Theorem~\ref{lem:algebraic-non-solvable}) is that a non-solvable finite group 
generated by an algebraically hyperelliptic collection is of a very special 
type, in particular $G$ has a quotient by a normal solvable subgroup which is 
isomorphic to the direct product of at most two simple groups.  

The next step is to cover the hyperelliptic subgroups by a bounded number of 
conjugacy classes of solvable subgroups (i.e., to find a collection of solvable 
subgroups such that each element of odd prime order in a hyperelliptic subgroup 
has a conjugate in one of these solvable subgroups); this concept of a solvable 
cover of a finite group is central for the proofs in the present paper since a 
bound on the number of elements of such a cover implies that there is a bound 
on the number of conjugacy classes of hyperelliptic subgroups. Using the 
classification of finite simple groups, we prove, in Proposition~\ref{simple}, 
that in any finite simple group the hyperelliptic subgroups can be covered, up
to conjugacy, by at most four solvable subgroups. This result, together with 
the characterisation of groups generated by algebraically hyperelliptic 
collections, gives directly the existence of a universal bound of fifty-seven,
much larger than the bound of fifteen obtained in Theorem~\ref{th:main4} by
exploiting extra topological considerations. Note that the existence of a 
universal (although non explicit) bound is ensured by the existence of only a 
finite number of simple sporadic groups.

In Section~\ref{sec:proof} we conclude the proof of Theorem~\ref{th:main4}. We 
can suppose that $G$ is generated by an algebraically hyperelliptic collection 
and that $G$ is not solvable. Under these hypotheses, the proof is divided into 
two cases. 

In the first case, we suppose that $G$ contains no rotation of order two. By  
geometric motivations, an involution acting dihedrally by conjugation on a 
hyperelliptic subgroup is a rotation. The absence of this type of involutions 
induces further restrictions on the structure of $G$; in particular, up to a 
quotient by a solvable subgroup, $G$ is a single simple group. By using these 
properties of $G$ and the solvable covers we prove that four is the upper bound 
in this case.

In the second case, $G$ contains a rotation of order two. The groups that may 
act on $3$-manifolds containing such an involution are listed in \cite{Mec}. 
We combine this with our result about groups generated by an algebraically  
hyperelliptic collection and we obtain that a quotient of $G$ by a solvable 
subgroup must be isomorphic to one of $\mathbb{A}_8$, $PSL_2(q)$ or 
$PSL_2(q)\times PSL_2(q')$. For these groups we explicitly find a solvable 
cover with a bounded number of elements and hence a universal bound as in 
Theorem~\ref{th:main4}.

%%%%%%%%%%%%%%%%%%%%%%%%%%%%%%%%%%%%%%%%%%%%%%%%%%%%%%%%%%%%%%%%%%

\section{Rotations and their properties}\label{sec:rotations}

In this section we shall establish some properties of rotations in general and
hyperelliptic ones in particular. 

\begin{remark}\label{rem:coverings}
Assume that $\psi$ is a hyperelliptic rotation acting on a $3$-manifold $M$.
The natural projection from $M$ to the space of orbits $M/\psi$ of $\psi$ 
is a cyclic cover of $\S^3$ branched along a knot $K=Fix(\psi)/\psi$. The 
converse is also true, that is any deck transformation generating the
automorphism group of a cyclic covering of $\S^3$ branched along a knot is a 
hyperelliptic rotation.
\end{remark}

We observe that cyclic branched covers of prime order are closely related to
$\QQ$-homology $3$-spheres.

\begin{remark}\label{rem:Qsphres}
\begin{enumerate}
\item[]
\item If the order of $\psi$ is a prime $p$, then $M$ is a $\ZZ/p$-homology
sphere \cite{Go}.
\item By Smith theory, if $f$ is a periodic diffeomorphism of order $p$, a
prime number, acting on a $\ZZ/p$-homology sphere, then $f$ either acts freely 
or is a rotation. 
\end{enumerate}
\end{remark}

We start with a somehow elementary remark which is however central to determine
constraints on finite groups acting on $3$-manifolds.

\begin{remark}\label{rem:normaliser}
Let $G\subset Diff^+(M)$ be a finite group of diffeomorphisms acting on a 
$3$-manifold $M$. One can choose a Riemannian metric on $M$ which is invariant 
by $G$ and with respect to which $G$ acts by isometries. Let now $\psi\in G$ be 
a rotation. Since the normaliser $\N_G(\langle \psi \rangle)$ of $\psi$ in $G$ 
consists precisely of those diffeomorphisms that leave the circle $Fix(\psi)$ 
invariant, we deduce that $\N_G(\langle \psi\rangle )$ is a finite subgroup of 
$\ZZ/2\ltimes(\QQ/\ZZ\oplus\QQ/\ZZ)$, where the nontrivial element in $\ZZ/2$ 
acts by conjugation sending each element of $\QQ/\ZZ\oplus\QQ/\ZZ$ to its 
inverse. Note that the elements of $\N_G(\langle \psi \rangle)$ are precisely 
those that rotate about $Fix(\psi)$, translate along $Fix(\psi)$, or invert the 
orientation of $Fix(\psi)$; in the last case the elements have order $2$ and 
non-empty fixed-point set meeting $Fix(\psi)$ in two points.

%Note that if $M\neq\S^3$ and $\psi$ is a hyperelliptic rotation of order $n>2$, 
%then its centraliser $\C_G(\langle \psi \rangle)$ in $G$ satisfies
%$1\longrightarrow \langle \psi \rangle \longrightarrow \C_G(\langle \psi
%\rangle) \longrightarrow H \longrightarrow 1$, where $H$ is cyclic, possibly 
%trivial. This follows easily from the positive solution to the Smith conjecture 
%which implies that any group of symmetries of a non-trivial knot $K$ (that is, 
%any finite group of orientation-preserving diffeomorphisms of $\S^3$ acting on 
%the pair $(\S^3,K)$) is either cyclic or dihedral. Moreover, since the 
%symmetries of a knot not acting freely have connected fixed-point set, the 
%possible elements of $\N_G(\langle \psi \rangle)\setminus 
%\C_G(\langle \psi \rangle) $ are rotations of order two.
\end{remark}

\begin{remark}\label{rem:smith}
Let us consider the $3$-sphere $\S^3$. According to Smith's theory, an
orientation-preserving finite-order diffeomorphism of $\S^3$ is a rotation if
and only if its fixed-point set is non-empty. Because of the positive solution 
to the Smith conjecture the fixed-point set of a rotation of $\S^3$ is the 
trivial knot. Morover, any group of symmetries of a non-trivial knot $K$ 
(that is, any finite group of orientation-preserving diffeomorphisms of $\S^3$ 
acting on the pair $(\S^3,K)$) is either cyclic or dihedral.  
\end{remark}

%\begin{definition}\label{def:rot-inv}
%With the notation of the above remark, we shall call 
%\emph{$Fix(\psi)$-rotations} the elements of $\N_G(\langle \psi \rangle)$ that 
%preserve the orientation of $Fix(\psi)$ and \emph{$Fix(\psi)$-inversions} those
%that reverse it.
%\end{definition}

\begin{Lemma}\label{lem:commuting powers suffice}
Let $\varphi$ and $\psi$ be two rotations contained in a finite group of 
orientation-preserving diffeomorphisms of a $3$-manifold $M$.
\begin{enumerate}
\item
A non-trivial power of $\psi$ of order different from $2$ commutes with a 
non-trivial power of $\varphi$, if and only if $\varphi$ and $\psi$ commute.
\item
Assume $M\neq\S^3$. If $\varphi$ and $\psi$ are hyperelliptic and 
$Fix(\varphi)=Fix(\psi)$, then $\langle \varphi \rangle=\langle\psi \rangle$ 
(in particular they have the same order).
\item
Assume $M\neq\S^3$. If $\varphi$ and $\psi$ are hyperelliptic, then 
$\langle \varphi \rangle$ and $\langle\psi \rangle$ are conjugate if and only 
if some non-trivial power of $\varphi$ is conjugate to some non-trivial power 
of $\psi$.
\end{enumerate}
\end{Lemma}

\demo

\noindent{\bf Part 1}
The sufficiency of the condition being obvious, we only need to prove the
necessity. Remark that we can assume that both rotations act as isometries for
some fixed Riemannian metric on the manifold. Denote by $f$ and $g$ the non
trivial powers of $\varphi$ and $\psi$, respectively. Note that, by 
definition, $Fix(\psi)=Fix(g)$ and $Fix(\varphi)=Fix(f)$. Since $g$ and $f$ 
commute, $g$ leaves invariant $Fix(\varphi)=Fix(f)$ and thus normalises every 
rotation about $Fix(\varphi)$. Moreover $g$ and $\varphi$ commute, for the 
order of $g$ is not $2$ (see Remark~\ref{rem:normaliser}). In particular, 
$\varphi$ leaves $Fix(\psi)=Fix(g)$ invariant and normalises every rotation 
about $Fix(\psi)$. The conclusion follows.

\noindent{\bf Part 2}
Reasoning as in Part 1, one sees that the two rotations commute. Assume, by
contradiction, that the subgroups they generate are different. Under this
assumption, at least one of the two subgroups is not contained in the other.
Without loss of generality we can assume that $\langle \varphi \rangle 
\not\subset \langle \psi \rangle$. Take the quotient of $M$ by the action of 
$\psi$. The second rotation $\varphi$ induces a non-trivial rotation of 
$\S^3$ which leaves the quotient knot $K=Fix(\psi)/\psi\subset\S^3$ invariant. 
Moreover, this induced rotation fixes pointwise the knot $K$. The positive 
solution to the Smith conjecture implies now that $K$ is the trivial knot and 
thus $M=\S^3$, against the hypothesis. 

\noindent{\bf Part 3} follows from 2 since the conjugate of a hyperelliptic
rotation is again a hyperelliptic rotation.
\qed

%\smallskip

%We notice also that a conjugate of a hyperelliptic rotation is a hyperelliptic 
%rotation. 

\begin{Corollary}\label{cor:powers suffice}
Let $G$ be a finite group of orientation-preserving diffeomorphisms acting on a
$3$-manifold $M\neq \S^3$. Let $\psi$ be a rotation and let $f\in G$ be an 
element of odd prime order which is a power of $\psi$. Then we have 
$\C_G(\langle f\rangle)=\C_G(\langle \psi\rangle)$ and $\N_G(\langle
f\rangle)=\N_G(\langle \psi\rangle)$.
\end{Corollary}

%The following two results show that the  presence of a hyperelliptic rotation 
%of order a multiple of a prime $p$ forces the Sylow $p$-subgroup to be of very 
%special type.

%\begin{definition}
%Let  $G \subset Diff^+(M)$ be a finite group acting on a closed orientable 
%$3$-manifold $M$. We say that an odd prime number $p$ is \emph{hyperelliptic} 
%for $G$ if  $p$ divides the order of a hyperilliptic rotation in $Diff^+(M).$
%\end{definition}

There is a natural bound on the number of hyperelliptic subgroups of order not 
a power of two which commute pairwise; we begin analysing the situation of the 
symmetry group of a knot.

\begin{definition}
A \emph{rotation of a knot} $K$ in $\S^3$ is a rotation $\psi$ of $\S^3$ such 
that $\psi(K)=K$ and $K\cap Fix(\psi)=\emptyset$. We shall say that $\psi$ is a 
\emph{full rotation} if $K/\psi$ in $\S^3=\S^3/\psi$ is the trivial knot.
\end{definition}

\begin{remark}\label{rem:full iff hyperelliptic}
Let $\psi$ and $\varphi$ be two commuting rotations acting on some manifold $M$ 
and with orders not both equal to $2$. 
%and whose axes are distinct. 
Assume that $\psi$ is hyperelliptic and $\varphi$ is not a power of $\psi$.
According to Remark~\ref{rem:normaliser} we have two situations: Either 
$Fix(\psi)\cap Fix(\varphi)=\emptyset$ and $\varphi$ induces a rotation $\phi$ 
of $K=Fix(\psi)/\psi$, or $Fix(\psi)=Fix(\varphi)$. In the former situation we 
have that $\varphi$ is hyperelliptic if and only if $\phi$ is a full rotation. 
This can be shown by considering the quotient of $M$ by the action of the group 
generated by $\psi$ and $\varphi$. This quotient is $\S^3$ and the projection 
onto it factors through $M/\varphi$, which can be seen as a cyclic cover of 
$\S^3$ branched along $K/\phi$. By the positive solution to the Smith 
conjecture, $M/\varphi$ is $\S^3$ if and only if $K/\phi$ is the trivial knot. 
In the latter situation, again because of the positive solution to Smith's 
conjecture, we have $M=\S^3$ (compare also Part 2 of 
Lemma~\ref{lem:commuting powers suffice}).
\end{remark}

The following finiteness result about commuting rotations of a non-trivial knot 
in $\S^3$ is one of the main ingredients in the proof of Theorem~\ref{th:main4} 
%and Theorem~\ref{th:main3} 
(see \cite[Proposition 2]{BoPa}, and \cite[Theorem 2]{BoPa} for a stronger 
result where commutativity is not required).

\begin{Proposition}\label{prop:commuting rotations} 
Let $K$ be a non-trivial knot in $\S^3$. Let us consider a set of pairwise 
commuting full rotations in $Diff^+(\S^3, K)$. The elements of the set 
generate at most two pairwise distinct cyclic subgroups. 
\end{Proposition}

\demo
Note, first of all, that, according to Remark~\ref{rem:smith}
%the positive solution to Smith's conjecture, 
the finite subgroup of $Diff^+(\S^3, K)$ generated by a finite set of pairwise 
commuting full rotations is cyclic.

Assume now, by contradiction, that there are three pairwise distinct cyclic 
subgroups generated by commuting full rotations of $K$, $\varphi$, $\psi$ and 
$\rho$ respectively. Note that such cyclic subgroups have distinct orders. 
Assume that two of them -say $\varphi$, $\psi$- have the same axis. 
%Note that, by hypothesis, they cannot have the same order. 
Fix the one with smaller order 
-say $\psi$-: since $\psi$ is a full rotation, the quotient $K/\psi$ is the trivial knot, and $\varphi$ induces 
a rotation of $K/\psi$ which is non-trivial since $\varphi$ commutes with 
$\psi$ and its order is larger than that of $\psi$. The axis $A$ of this 
induced symmetry is the image of $Fix(\psi)$ in the quotient $\S^3/\psi$ by the 
action of $\psi$. In particular $K/\psi$ and $A$ form a Hopf link and $K$ is 
the trivial knot: this follows from the equivariant Dehn lemma, see \cite{Hil}. 

We can thus assume that the axes are pairwise disjoint. Indeed, this follows
from Remark~\ref{rem:normaliser} taking into account that the rotations commute
pairwise and at most one of them has order $2$.
%Since the rotations commute, even if one of them has order $2$, it cannot act 
%as a strong inversion on the axes of the other rotations. 
Therefore we would have that the axis of $\rho$, which is a trivial knot, 
admits two commuting rotations, $\varphi$ and $\psi$, with distinct axes, which 
is impossible: this follows, for instance, from the fact (see 
\cite[Thm 5.2]{EL}) that one can find a fibration of the complement of the 
trivial knot which is equivariant with respect to the two symmetries.
\qed
\smallskip

Observe that the proof of the proposition shows that two commuting full
rotations of a non-trivial knot either generate the same cyclic subgroup or 
have disjoint axes. 
%This follows also from the analogous conclusion holding for
%hyperelliptic rotations.

%\begin{remark}\label{rem:prime} 
We stress that if a knot $K \subset \S^3$ admits a full rotation, then it is a 
prime knot, see \cite[Lemma 2]{BoPa}.
%\end{remark}

\begin{Lemma}\label{lem:abelian-case}
Let $G$ be a finite group of diffeomorphisms acting on a $3$-manifold
$M\neq\S^3$ and $\{H_1,\dots,H_m\}$ be a set of hyperelliptic subgroups of $G$ 
of order not a power of two. Suppose that there exists an abelian subgroup of 
$G$, containing at least an element of odd prime order of each $H_i$, then 
$m\leq 3$. Moreover either the orders of the $H_i$ are pairwise coprime 
%distinct 
or $m\leq 2.$
\end{Lemma} 

\demo 

By Lemma~\ref{lem:commuting powers suffice} we obtain that the subgroups $H_i$ 
commute and have trivial intersection pairwise. Consider the cyclic branched 
covering $M\rightarrow M/H_1\cong \S^3$ over the knot $Fix(H_1)/H_1$. By 
projecting $H_i$ with $i\geq 2$ to $M/H_1$ we obtain full rotations of 
$Fix(H_1)/H_1.$ Note that if two subgroups $H_i$ have the same order that is
moreover a divisor of the order of $H_1$, a priori they might map to the same 
subgroup in $\langle H_1,\dots,H_m \rangle/H_1$. We claim, however, that the 
induced full rotations are distinct so there are $m-1$ of them. This follows 
from the fact that, for different indices $i$ and $j$, $Fix(H_i)$ and 
$Fix(H_j)$ are disjoint according to Part 2 of 
Lemma~\ref{lem:commuting powers suffice}. Now, since the subgroups commute, for 
all $i=1,\dots, m$ we have $H_1(Fix(H_i))=Fix(H_i)$, so that the fixed-point 
sets of the induced full rotations are disjoint, too, and the full rotations 
are pairwise distinct. By Proposition~\ref{prop:commuting rotations} we obtain 
$m-1 \leq 2$. Note that, by Remark~\ref{rem:smith}
%the positive solution to the Smith conjecture, 
a non-trivial knot cannot admit two distinct and commuting cyclic 
groups of symmetries of the same order. This proves the latter part of the 
lemma.
\qed

\smallskip

The above lemma implies directly the following corollary.

\begin{Corollary}\label{cor:2hyprot}
Let $p$ be an odd prime and assume that $H\cong \ZZ/p\oplus\ZZ/p$ acts on a 
$3$-manifold $M\neq\S^3$, the group $H$ contains at most two distinct cyclic 
subgroups generated by non-trivial powers of hyperelliptic rotations. 
\end{Corollary}

%Note that this is not immediately obvious for $\ZZ/p\oplus\ZZ/p$ contains $p+1$ 
%cyclic subgroups with pairwise trivial intersection and $p$ of them project to 
%the same cyclic subgroup in the quotient by the remaining one.

%Indeed by a slighty different proof we could obtain a more general result, in fact if the group $H$ contains  the power of a hyperelliptic  involution,  it contains  at most two distinct  cyclic subgroups generated by  a rotation of order $p$.

We now collect several facts relative to hyperelliptic rotations that can be
deduced from the discussion in this section.

\begin{remark}\label{rem:all about hr}
Let $G$ be a finite group of orientation-preserving diffeomorphisms of a
$3$-manifold $M\neq\S^3$. Let $\psi\in G$ be a hyperelliptic rotation of order 
not a power of $2$. Recall that the structures of the centraliser and of the 
normaliser of $\langle \psi \rangle$ are described in 
Remark~\ref{rem:normaliser}. Let $N$ denote the normaliser 
$N_G(\langle \psi \rangle)$.

\begin{enumerate}
\item The centraliser $\C_G(\langle \psi \rangle)$ of $\psi$ in $G$ satisfies
$1\longrightarrow \langle \psi \rangle \longrightarrow \C_G(\langle \psi
\rangle) \longrightarrow H \longrightarrow 1$, where $H$ is cyclic, possibly
trivial. This follows from Remark~\ref{rem:smith}. 
\item Since the symmetries of a knot not acting freely have connected 
fixed-point set (see again Remark~\ref{rem:smith}), the possible elements of 
$N\setminus \C_G(\langle \psi \rangle) $ are rotations of order two.
\item According to last part of Lemma~\ref{lem:abelian-case}, the centraliser 
$\C_G(\langle \psi \rangle)$ contains at most one more cyclic subgroup of the 
same order as $\langle \psi \rangle$ and generated by a hyperelliptic rotation. 
\item Observe that we have $\N_G(N)=\N_G(\C_G(\langle \psi\rangle))$. 
\item Of course, the conjugate of a hyperelliptic rotation is again a
hyperelliptic rotation. Consider an element $g\in \N_G(N)$. Then either $g$ 
conjugates $\langle \psi \rangle$ to itself and so $g$ belongs to $N$, or $g$ 
does not normalise $\langle \psi \rangle$. In this case, 
$\C_G(\langle \psi\rangle)$ contains precisely two cyclic subgroups generated 
by hyperelliptic rotations of the same order as $\psi$ ($\langle \psi\rangle$ 
and another one) and $g$ exchanges them. It follows that 
$\C_G(\langle \psi \rangle)$ has index a divisor of $4$ in $\N_G(N)$.
\item Let $f\in G$ be an element of odd prime order which is a power of 
$\psi$. Then we have $\C_G(\langle f\rangle)=\C_G(\langle \psi\rangle)$ and 
$\N_G(\langle f\rangle)=N$. This is just Corollary~\ref{cor:powers suffice}.
\end{enumerate}
\end{remark}

%\begin{remark}\label{rem:normnorm}
%Let $G$ be a finite group of diffeomorphism acting on a manifold $M\neq\S^3$. 
%Let $\psi\in G$ be a hyperelliptic rotation of order not a power of two and let 
%$N$ be the normaliser of $\langle \psi \rangle$ in $G$. It follows from 
%Remark~\ref{rem:normaliser} and Corollary~\ref{cor:2hyprot} that $N$ has index 
%at most two in $\N_G(N)$; moreover the index is precisely two if and only if 
%$N$ contains two distinct cyclic subgroups conjugate to $\langle \psi \rangle$ 
%and every element in $\N_G(N)\setminus N$ acts by exchanging the two subgroups 
%(see also Proposition~\ref{pro:hyperellipticnormsylow} together with
%Lemma~\ref{lem:commuting powers suffice}: observe that
%$\N_G(N)=\N_G(\C_G(\langle \psi\rangle))$). 
%\end{remark}

The following general observation will be useful in the sequel.

\begin{remark}\label{rem:conjugate-rotations}
With the notation used in the previous remark, assume that there are elements 
in $\N_G(N)$ which do not act in the same way on $\C_G(\langle \psi\rangle)$, 
then necessarily the index of $N$ in $\N_G(N)$ is two and 
$\C_G(\langle \psi\rangle)$ contains two distinct subgroups conjugate to 
$\langle \psi\rangle$.
\end{remark}

We end this section by providing a dictionary translating between algebraic
properties of the structure of $\N_G(N)$ and the symmetries of the knot
$K=Fix(\psi)/\psi$. This will not be needed in the proofs of our results but
will provide a geometric interpretation of the different situations that occur
in the proof of Theorem~\ref{th:main4} (see Remark~\ref{r:casesbysymmetry}). 
We start with the following definition.

\begin{definition}
A $2$-component link is called \emph{exchangeable} if there exists an
orientation-preserving diffeomorphism of $\S^3$ which exchanges the
two components of the link.

Let $K$ be a knot and $\rho$ a rotation of $K$ of order $n$ and with axis $A$.
Consider the $2$-component link $\overline{K} \cup \overline{A}$ consisting of
the images of the knot $K$ and of the axis $A$ in the quotient
$\S^3/ \rho$ of the $3$-sphere by the action of $\rho$. Note
that at least one component of this link (i.e. $\overline{A}$) is trivial. We
call $K$ \emph{$n$-self-symmetric} if $\overline{K} \cup \overline{A}$ is
exchangeable. In this case $\rho$ is a full rotation of $K$.
\end{definition}

\begin{figure}[H]
\begin{center}
 {
  \includegraphics[height=3.5cm]{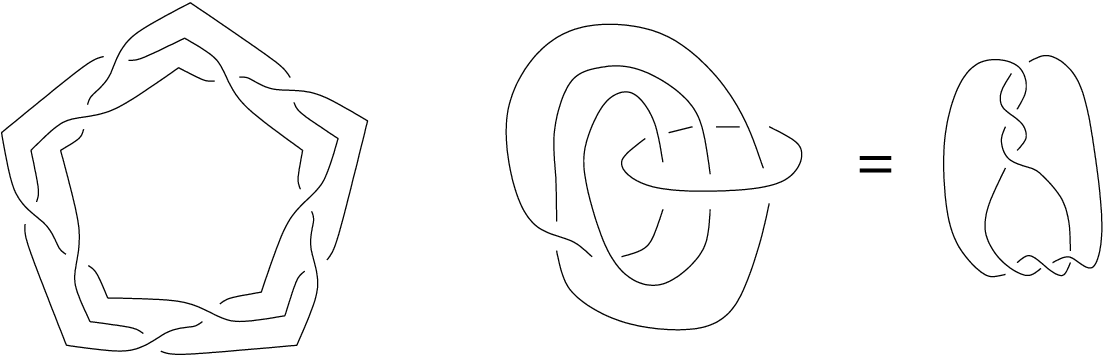}
 }
\end{center}
\caption{A $5$-self-symmetric knot on the left, and its exchangeable quotient
link on the right.}
\label{fig:selfsymmetric}
\end{figure}

Since the structure of the normaliser of $\langle \psi \rangle$ and of its
centraliser, only depend on the symmetries of $K$ that lift to $G$, we 
introduce the following definitions:

\begin{definition}\label{def:relG}
Let $G$ be a finite group of orientation preserving diffeomorphisms of a
closed connected $3$-manifold $M$. Let $\psi$ be a hyperelliptic rotation
contained in $G$ with quotient knot $K$. We say that $K$ is \emph{strongly
invertible with respect to $G$} if $K$ admits a strong inversion that lifts to
$G$. Similarly we say that $K$ is \emph{self-symmetric with respect to $G$} if
$G$ contains an element $\psi'$ conjugate to $\psi$ such that the
subgroup $\langle \psi,\psi'\rangle$ is abelian of rank $2$, i.e. not
cyclic. Remark that in the latter case $K$ is $n$-self-symmetric, where $n$ is
the order of $\psi$.
\end{definition}

The proof of the following facts is elementary and left to the reader.

\begin{Proposition}\label{prop:symmetries and normalizer}
Let $K$, $\psi$, $N$, and $G$ be as above.
\begin{itemize}
\item The centraliser of $\langle \psi \rangle$ in $G$ is contained with index 
$2$ in $N$ if and only if $K$ is strongly invertible with respect to $G$.
\item $N$ is contained with index $2$ in $\N_G(N)$ if and only if $K$ is 
self-symmetric with respect to $G$.
\end{itemize}
Moreover, if $M$ is hyperbolic and $G=Iso^+(M)$, $K$ is strongly invertible if
and only if it is strongly invertible with respect to $G$, and it is
self-symmetric with respect to $G$ if and only if it is $n$-self-symmetric,
where $n$ is the order of the hyperelliptic rotation $\psi$.
\end{Proposition}

%\bigskip

%%%%%%%%%%%%%%%%%%%%%%%%%%%%%%%%%%%%%%%%%%%%%%%%%%%%%%%%%%%%%%%%%%%%%%%%%%%

\section{Group theoretical results}\label{sec: algebraic lemmas}

%In this section we provide some algebraic Lemmas which we need in the proof of 
%Theorem~\ref{th:main4}.

This section is devoted to the proofs of the group-theoretical results that are 
used to obtain the bound provided by Theorem~\ref{th:main4}.

The first result, Theorem~\ref{lem:algebraic-non-solvable}, is probably the
more interesting of the two from a group-theoretical point of view. It
describes the finite groups $G$ that may be generated by hyperelliptic 
rotations. In this settings, a hyperelliptic rotation is simply an element of
$G$ of order not a power of two, satisfying a purely algebraic condition (see 
Definition~\ref{ah}, below) on the structure of the normaliser of (a power of
odd prime order of) the element. More precisely, 
Theorem~\ref{lem:algebraic-non-solvable} states that either $G$ is solvable, or
it admits a quotient by a normal solvable subgroup so that the quotient is
either the product of a simple group times a cyclic group of odd order, or the 
product of two simple groups. We see that the fact of being generated by
hyperelliptic rotations puts very strong constraints on the structure of $G$.

\medskip

The second result, Proposition~\ref{simple}, although possibly not as striking 
as the previous one, is however key to be able to carry out the strategy of 
bounding the number of hyperelliptic subgroups by covering them with solvable
subgroups (see Definition~\ref{def:solcov}). Indeed, in 
Proposition~\ref{simple} we establish the fact that four conjugacy classes of 
solvable subgroups are sufficient to contain, up to conjugacy, all 
hyperelliptic subgroups of orders not a power of two that may sit inside a 
finite simple group.

Note that these two main results imply that $4\times 4$ solvable subgroups 
suffice to cover all hyperelliptic groups in $G$, and together with 
Lemma~\ref{lem:algebraic-solvable}, and Lemma~\ref{lem:abelian-case} show that 
a finite group $G$ can contain at most forty-eight conjugacy classes of 
hyperelliptic subgroups of orders that are not powers of $2$. 

\bigskip

In the following $G$ will denote a finite group. We start with some preliminary
observations and definitions.

%\begin{definition}
%Let $S=\{H_1,\,H_2,\dots,H_u\}$ be a collection of subgroups of $G$, the collection $S$ \textit{invariably generates} $G$ 
%if for each $g_1,\,g_2,\dots g_u \in G$ the subgroups  $\{H_1^{g_1}\,H_2^{g_2},\dots,H_u^{g_u}\}$ generate $G.$
%\end{definition}

\begin{definition}\label{ah}
A collection $\{C_i\}$ of subgroups of $G,$ each of odd prime order $p_i$, is said to be 
\emph{algebraically hyperelliptic} if, for each $i$, the following conditions are satisfied:
\begin{enumerate} 
\item the centraliser of $C_i$ in $G$ is abelian of rank at most two and has index at 
most two in the normaliser of $C_i$ in $G$;
\item each element belonging to the normaliser of $C_i$  but not to the 
centraliser inverts by conjugation each element in the centraliser;
\item if $C_i$ is contained in a Sylow $p_i$-subgroup $S_i$ of $G$, then $S_i$ 
contains at most two distinct conjugates of $C_i$.
\end{enumerate}
\end{definition}

We remark that in this definition the primes $p_i$ are not necessarily 
pairwise distinct.

\begin{Proposition}\label{pro:hyperellipticnormsylow}
Let $S$ be a Sylow $p_i$-subgroup where $p_i$ is the order of a group $C_i$ 
belonging to an algebraically hyperelliptic collection. 
\begin{enumerate}
\item $S$ is either cyclic or the product of two cyclic groups, and 
\item $\N_G(S)$ contains with index at most $2$ the normaliser of a conjugate 
of $C_i$, and contains an abelian subgroup of rank at most $2$ with index a 
divisor of $4$. In particular $\N_G(S)$ is solvable.
\end{enumerate}
\end{Proposition}

\demo

Up to conjugation we can suppose that $S$ contains $C_i$. By  Properties~1 and 2 
in Definition~\ref{ah} and the fact that $p_i$ is odd, the normaliser $\N_S(C_i)$ is abelian of rank at most 
two. Property~3 in Definition~\ref{ah} implies that 
$\N_S(\N_S(C_i))=\N_S(C_i)$. Since $S$ is a $p_i$-subgroup we obtain that 
$S=\N_S(C_i)$ and we get the first part of the thesis. 

Since the subgroup $\N_G(S)$ normalises the maximal elementary abelian subgroup 
of $S$, we obtain also the second part of the thesis.
\qed

\bigskip

%The following general observation will be useful in the sequel.

%\begin{remark}\label{rem:conjugate-rotations}
%We work under the hypotheses of Proposition~\ref{pro:hyperellipticnormsylow} 
%and we suppose that $C_i$ is contained in the Sylow $p_i$-subgroup $S$. The 
%normaliser of $S$ contains the normaliser of $C_i$ with index $2$ if and only 
%there exist an element $g\in N_G(S)$ such that $C_i^g\neq C_i.$  This case 
%happens if and only if $\N_G(S)$ contains elements of order a power of $2$ 
%which do not act in the same way on all elements of order $p_i$ in $S$. Indeed, 
%all elements in $\N_G(C_i)$ either commute with all elements of order $p_i$ or 
%act dihedrally. On the other hand, any element $h$ in 
%$\N_G(S)\setminus \N_G(C_i)$ conjugates $C_i$ to $C_i^g.$ If $f$ is a 
%genarator of $C_i$ and $f'=f^h$ (this element is a generator of $C_i^g$), we 
%obtain that $hf'fh^{-1}=f'f$ and $hf'f^{-1}h^{-1}=(f'f^{-1})^{-1}$, i.e. $h$ 
%acts dihedrally on some elements of order $p_i$ while it commutes with others.
%\end{remark}

\begin{Lemma}\label{lem:algebraic-solvable}
Let $G$ be a solvable group containing an algebraically hyperelliptic 
collection $\{C_1,\cdots , C_m\}$ of subgroups of odd prime order. Then there 
exists an abelian subgroup of $G$ containing a conjugate of $C_i$, for each 
$i=1,\dots,m.$ In particular the subgroups $C_i$ commute pairwise up to 
conjugacy. 
\end{Lemma}

\demo
Let $\pi$ be the set of the orders of the $C_i$ and let $B$ be a Hall 
$\pi$-subgroup of $G$. Each $C_i$ is conjugate to a subgroup of $B$. Since 
$\pi$ contains only odd primes, Definition~\ref{ah} yields that centraliser and 
normaliser of every Sylow $p$-subgroup of $B$ coincide. By Burnside's 
$p$-complement theorem (see~\cite[Theorem 2.10 page 144]{Su}), every Sylow 
$p$-subgroup of $B$ has a normal complement, and hence $B$ is abelian. 
\qed

\begin{remark}\label{rem:normalizer-coprime-quotient}
Suppose that $N$ is a normal subgroup of $G$ and $H$ is a $p$-subgroup of $G$. 
If the order of $N$ is coprime with $p$, then the normaliser of the projection 
of $H$ to $G/N$ is the projection of the normaliser of $H$ in $G$, that is 
$$\N_{G/N}(HN/N)=\N_G(H)N/N.$$
The inclusion $\supseteq$ holds trivially. We prove briefly the other 
inclusion. Let $fN$ be an element of $\N_{G/N}(HN/N)$, then $H^f\subseteq HN$. 
Both $H^f$ and $H$ are Sylow  $p$-subgroups of $HN$ and by the second Sylow 
theorem they are conjugate by an element $hn\in HN.$ We obtain that 
$H^{fhn}=H$, and hence $f\in \N_G(H)N.$

Analogously if $f$ is an element of prime order coprime with the order of $N$, 
then $\C_{G/N}(fN)=\C_G(f)N/N.$
\end{remark}

Recall that a finite group $Q$ is \emph{quasisimple} if it is perfect (the
abelianised group is trivial) and the factor group $Q/\Z$ of $Q$ by its centre
$\Z$ is a nonabelian simple group (see \cite[chapter 6.6]{Su}). A group $E$ is 
\emph{semisimple} if it is perfect and the factor group $E/\Z(E)$ is a direct 
product of nonabelian simple groups. A semisimple group $E$ is a central 
product of quasisimple groups which are uniquely determined. Any finite group 
$G$ has a unique maximal semisimple normal subgroup $E(G)$ (maybe trivial), 
which is characteristic in $G$. The subgroup $E(G)$ is called the \emph{layer} 
of $G$ and the quasisimple factors of $E(G)$ are called the \emph{components} 
of $G$.

The maximal normal nilpotent subgroup of a finite group $G$ is called the 
\emph{Fitting subgroup} and is usually denoted by $F(G)$. The Fitting subgroup 
commutes elementwise with the layer of $G$. The normal subgroup generated by 
$E(G)$ and by $F(G)$ is called the \emph{generalised Fitting subgroup} and is
usually denoted by $F^*(G)$. The generalised Fitting subgroup has the 
important property to contain its centraliser in $G$, which thus coincides 
with the centre of $F^*(G)$. For further properties of the generalised Fitting 
subgroup see \cite[Section 6.6.]{Su}.

\begin{Theorem}\label{lem:algebraic-non-solvable}
Suppose that $G$ is generated by the algebraically hyperelliptic collection 
$\mathcal H:=\{C_1,\cdots,C_m\}$. Denote by $p_i$ the order of $C_i$ and by $A$ 
the maximal normal solvable subgroup of order coprime with every $p_i$.   

If $G$ is non-solvable, then the following properties hold:

\begin{enumerate}

\item  every  $p_i$ divides the order of any component of $G/A;$

\item either $G/A$ is the direct product of a cyclic group of odd order and a 
simple group or it is the direct product of two simple groups;

\item if in addition $G$ does not contain any involution acting dihedrally on 
any $C_i$, then $E(G/A)$ is simple, every $p_i$ divides the order of $F(G/A)$ 
and a Sylow $p_i$-subgroup of $G$ contains exactly two conjugates of $C_i$.

\end{enumerate}

\end{Theorem}
%[Remark: in this case we can assume  generation instead invariant generation]

\demo
Let $\pi$ be the set of the primes $p_i.$ 

By Remark~\ref{rem:normalizer-coprime-quotient} we can suppose that $A$ is 
trivial and $F(G)$ is a $\pi$-group.

\smallskip
\begin{Claim}\label{c1}
$F(G)$ is cyclic and $E(G)$ is not trivial.
\end{Claim}

Suppose first that $F(G)$ contains an abelian $p_i$-subgroup $S$ of rank two.  
Then $S$, being the maximal elementary abelian $p_i$-subgroup contained in 
$F(G)$, is normal in $G$ and contains $C_i$. This implies that $G$ is solvable 
and we get a contradiction. Hence $F(G)$ is cyclic.

If $E(G)$ is trivial, then $F(G)=F^*(G)$ and $G/F(G)$ is isomorphic to a 
subgroup of Aut$F(G).$ Since $F(G)$ is cyclic Aut$F(G)$ is solvable, and we get 
again a contradiction. 

\smallskip

\begin{Claim}\label{c2}
Each $p_i$ divides the order of any component of $G$. Moreover the components 
of $G$ are simple groups and are at most two.
\end{Claim}

Since the Sylow $p_i$-subgroups are abelian and $A$ is trivial, 
by~\cite[Exercise 1, page 161]{Su} the components of $G$ have trivial centre.

Now we prove that each component of $G$ is normalised by any $C_i$. Let $f_i$ 
be a generator of $C_i$ and $Q$ a component of $G.$ Suppose by contradiction 
that $Q$ is not normalised by $f_i$. We define the following subgroup:
$$Q_c=\{x f_i x  f_i^{-1}\dots  f_i^{p_i-1} x  f_i^{p_i-1} \mid x \in Q\}.$$
Since the components of $G$ commute elementwise, $Q_c$ is a subgroup of $G$ 
isomorphic to $Q$. Moreover, each element of $Q_c$ commutes with $f_i$ and this 
gives a contradiction.

We have that $C_i$ normalises $Q$ but cannot centralise it, so the action by 
conjugation of $f_i$ on $Q$ is not trivial.

Assume that $Q$ is either sporadic or alternating. Since the order of the 
outer automorphism group of any such simple group is a (possibly trivial) 
power of $2$ (see \cite[Section 5.2 and 5.3]{GLS}), we conclude that $f$ must 
induce an inner automorphism of $Q$. In particular $p_i$ divides the order of 
$Q$.

We can thus assume that $Q$ is a simple group of Lie type.

Recall that, by \cite[Theorem~2.5.12]{GLS}, $Aut(Q)$ is the semidirect product 
of a normal subgroup $Inndiag(Q)$, containing the subgroup $Inn(Q)$ of inner 
automorphisms, and a group $\Phi\Gamma$, where, roughly speaking, $\Phi$ is the 
group of automorphisms of $Q$ induced by the automorphisms of the defining 
field and $\Gamma$ is the group of automorphisms of $Q$ induced by the 
symmetries of the Dynkin diagram associated to $Q$ (see \cite{GLS} for the 
exact definition). By \cite[Theorem~2.5.12.(c)]{GLS}, every prime divisor of 
$|Inndiag(Q)|$ divides $|Q|$. Thus we can assume that the automorphism 
induced by $f_i$ on $Q$ is not contained in $Inndiag(Q)$ and its projection 
$\theta$ on $Aut(Q)/Inndiag(Q)\cong \Phi\Gamma$ has order $p_i$. We will find a
contradiction showing that in this case the centraliser of $f_i$ in $Q$ is not 
abelian.  

Write $\theta=\phi\gamma$, with $\phi \in \Phi$ and $\gamma\in \Gamma$. If 
$\phi=1$, then $\gamma$ is nontrivial and $f_i$ induces a graph automorphism 
according to \cite[Definition 2.5.13]{GLS}. Since $p_i$ is odd, the only 
possibility is that $Q$ is $D_4(q)$ and $p_i=3$ (see 
\cite[Theorem 2.5.12 (e)]{GLS}). The centraliser of $f_i$ in $Q$ is nonabelian  
by \cite[Table~4.7.3 and Proposition 4.9.2.]{GLS}. If $\phi\neq 1$ and $Q$ is 
not isomorphic to the group $^3D_4(q)$, then the structure of the centraliser 
of $f_i$ in $Q$ is described by \cite[Theorem~4.9.1]{GLS}, and it is 
nonabelian. Finally, if $\phi\neq 1$ and $Q\cong {^3D_4(q)}$, the structure of 
the non abelian centraliser of $f_i$ in $Q$ follows from 
\cite[Proposition~4.2.4]{GLS}. We proved that the automorphism induced by $f_i$ 
is contained in $Inndiag(Q)$ and $p_i$ divides $|Q|$.

Since $p_i$ divides the order of any component, $G$ has at most two components.

\smallskip

\begin{Claim}\label{c3}
$G=E(G)F(G)$.
\end{Claim}

We prove first that $G=E(G) \C_G(E(G)).$ Let us assume by contradiction that 
there exists $C_i$ with trivial intersection with $E(G)\C_G(E(G))$ and denote 
by $f$ a generator of $C_i.$ Since $p_i$ divides the order of every component 
of $G$ and the Sylow $p_i$-subgroup has rank at most $2$, we get that $E(G)$ 
has only one component which we denote by $Q$. The Sylow $p_i$-subgroups of $Q$ 
are cyclic. Moreover, by the first part of the proof, the automorphism induced 
by $f$ on $Q$ is inner-diagonal. If it is inner, we obtain $f$ as a product of 
an element that centralises $Q$ and an element in $Q$, a contradiction to our 
assumption; otherwise, we get again contradiction, since, by 
\cite[Theorem~2.5.12]{GLS} and \cite[(33.14)]{A}, a group of Lie type with 
cyclic Sylow $p_i$-subgroup cannot have a diagonal automorphism of order $p_i$.

Hence, all the subgroups $C_i$ are contained in $E(G)\C_G(E(G))$ and, since 
they generate $G$, we obtain that $G=E(G)\C_G(E(G)).$ 

Now if $F(G)=1$, then $F^\ast(G)=E(G)$ and hence 
$\C_G(E(G))=\C_G(F^\ast (G))=\Z(E(G))=1$ and the claim is proved. So suppose 
that $F(G)\neq 1$. Then, since $F(G)$ is a $\pi$-group, there is at least one 
subgroup $C_i$ that is contained in $E(G)F(G)$. Hence there is a subgroup $T_1$ 
of $E(G)$ with order $p_i$ and a subgroup $T_2$ of $F(G)$ with order $p_i$ such 
that $C_i\leq T_1T_2$. Since $F(G)$ is cyclic, $T_2$ is normal in $G$ and so 
$\C_G(E(G))\leq \N_G(T_1T_2)$. But $\N_G(T_1T_2)$ has an abelian normal 
subgroup of index a divisor of $4$, so $\C_G(E(G))$ is solvable. This shows 
that every Sylow $p$-subgroup of $\C_G(E(G))$ for $p$ odd is contained in 
$F(G)$, whence it follows that $G=E(G)F(G)$.

%Since the centraliser of the generalised Fitting subgroup coincides with its 
%centre $\Z(F^*(G))$ and since $C_G(E(G))$ acts trivially on $E(G)$ by definition, 
%the quotient $C_G(E(G))/F(G)$ merges injectively in the automorphism group of 
%$F(G)$. Since the  Fitting subgroup is cyclic, then its automorphism group is abelian and $C_G(E(G))$ is solvable. This implies that $C_G(E(G))$ is a $\pi$-group ($A$ is trivial) and the Sylow subgroups of $C_G(E(G))$ are cyclic.

%Finally we prove that $C_G(E(G))$ is abelian, and hence we obtain $C_G(E(G))=F(G)$

%If $E(G)$ contains a Sylow $p$-subgroup,  $C_G(E(G))$ has to be abelian and we are done.

%Otherwise $E(G)$ has only one component $Q$ and the Sylow $p_i$-subgroups of $Q$ are non-trivial and cyclic. T
%By the solvability of $C_G(E(G))$  we have a subnormal series $\{1\}=H_O\triangleleft H_1 \triangleleft H_2\triangleleft\cdots \triangleleft H_a=C_G(E(G))$ such that $H_i/H_{i-1}$ is abelian. We suppose that $H_1$ is the maximal normal abelian subgroup of $H_2.$
%If $H_1=C_G(E(G))$ we are done, otherwise we find in $H_2$  an element $f$ of order $p_i ^\alpha$ that by conjugation acts non-trivially on $H_1$. In $H_1$ the exists an element of order $p_j$ with $p_j\neq p_i$ that is normalized by $f$ but does not commute with it. We obtain that $f$ normalize an elementary abelian $p_j$-subgroup of rank two but is not in its centralizer, and this is impossible.

\begin{Claim}\label{c4}
If no involution of $G$ acts dihedrally on any $C_i$, then $E(G)$ is simple, 
and, for every $i$, $p_i$ divides the order of $F(G)$ and a Sylow 
$p_i$-subgroup of $G$ contains exactly two conjugates of $C_i$.
\end{Claim}

Suppose no involution of $G$ acts dihedrally on any $C_i$. By 
Definition~\ref{ah}, it follows that $\N_G(C_i)=\C_G(C_i)$ for every $i=1, 
\ldots, m$. Let $Q$ be a component of $E(G)$ and suppose by contradiction that 
$Q$ contains a Sylow $p_i$-subgroup $S$ of rank two, for some $i$. Up to 
conjugation we can suppose that $S$ contains $C_i.$ By Definition~\ref{ah}, 
$\N_Q(S)$ contains with index at most two the abelian group $\C_Q(C_i)$. Since 
$Q$ is simple, Burnside's $p$-complement theorem (see~\cite[Theorem 2.10 page 
144]{Su}) yields that $\N_Q(S)$ is not abelian. Therefore, $\N_Q(S)$ contains 
with index two $\N_Q(C_i)$ and the elements of $\N_Q(S)\setminus \N_Q(C_i)$  
conjugate $C_i$ to a cyclic subgroup distinct from $C_i$. By using 
\cite[Exercise 1, page 161]{Su} and the fact that $Q$ is perfect we get again a 
contradiction. Hence, for every $i\in \{1, \ldots , m\}$, the Sylow 
$p_i$-subgroups of $Q$ are cyclic. Now, as above for every Sylow subgroup $S$ 
of $Q$ we have $\N_Q(S)\neq\C_Q(S)$. Since $\N_G(C_i)=\C_G(C_i)$ for every 
$i=1, \ldots, m$, that implies that the Sylow $p_i$-subgroups of $G$ are not 
cyclic and hence $p_i$ divides the order of $F(G)$ for every $i$. 
\qed
\smallskip

To bound above the number of conjugacy classes of hyperelliptic rotations, our 
strategy will consist of conjugating hyperelliptic rotations in solvable 
subgroups of $G$, where they are forced to commute, hence, in analogy with the standard definition of normal covers, we introduce the 
notion of solvable normal $\pi$-cover and we prove the following lemma.

\begin{definition}\label{def:solcov}
Let $G$ be a finite group. Let $\pi$ be a set of primes dividing $|G|$. We will 
call a collection $\C$ of subgroups of $G$ a \emph{solvable normal $\pi$-cover} 
of $G$ if every element of $G$ of prime order $p$ belonging to $\pi$ is 
contained in an element of $\C$ and for every $g\in G$, $H\in \C$ we have that 
$H^g\in \C$. We denote by $\gamma^s_{\pi}(G)$ the smallest number of conjugacy 
classes of subgroups in a solvable normal $\pi$-cover of $G$. Here, the letter ``$s$'' stands for ``solvable'', and is used to distinguish this number  from  $\gamma(G)$, that is the standard notation in the case of covers by subgroups that are not requested to be solvable.  Note that, since
Sylow subgroups are clearly solvable, $\gamma^s_{\pi}(G)\leq |\pi|$.
\end{definition}

\begin{Proposition}\label{simple}
Let $G$ be a finite nonabelian simple group. If $\pi$ is the set of odd primes 
$p$ such that $G$ has cyclic Sylow $p$-subgroup, the centraliser of $C_G(g)$ is 
abelian for every element $g\in G$ of order $p$ and the normaliser of any 
subgroup of order $p$ contains with index at most two its centraliser, then 
$\gamma^s_{\pi}(G)\leq 4$.
\end{Proposition}

\demo

If $G$ is a sporadic group, the primes dividing the order of the group do not 
satisfy the condition on the normaliser.

If $G$ is isomorphic to the alternating group $\mathbb{A}_n$ and $p\in \pi$, 
then the condition on the centraliser of the elements of order $p$ implies that $p>n-4$ and $\gamma^s_{\pi}(G)\leq |\pi|\leq 2.$

The only remaining case is that of groups of Lie type.

Let $G\cong {\Sigma}_n(q)$ or $G\cong {^d\Sigma}_n(q)$, where $q$ is a power of 
a prime $t$. Here we use the same notation as in \cite{GLS}: the symbol 
$\Sigma(q)$ (resp. $^d\Sigma(q)$) may refer to finite groups in different 
isomorphism classes, each of them is an untwisted (resp. twisted) finite group 
of Lie type with root system $\Sigma$ (see \cite[Remark 2.2.5]{GLS}). Any 
finite group of Lie type is quasisimple with the exception of the following 
groups: $A_1(2),\,A_1(3),\,^2A_2(2),\,^2B_2(2),\,B_2(2),\,G_2(2),\,^2F_4(2)$ 
and $^2G_2(3)$ (see \cite[Theorem 2.2.7]{GLS}).  

If $t\in \pi$, then by~\cite[Theorem~3.3.3]{GLS}, either $t=3$ and 
$G\cong (^2G_2(3))'$ or $G\cong A_1(s)$. In the former case the order of 
$^2G_2(3)'$ is divided only by two odd primes, thus $\gamma^s_{\pi}(G)\leq 2$; 
in the latter case we have $\gamma^s_{\pi}(G)\leq 2$ (see for example 
\cite{H}).

Assume now that $t\not\in \pi$. By \cite[Paragraph~4.10]{GLS}, since the 
Sylow subgroups are cyclic, every element of order $p\in\pi$ is contained in a 
maximal torus of $G$, and clearly a maximal torus is abelian. 

Therefore, we need only to bound the number of conjugacy classes of cyclic 
maximal tori in $G$ with abelian centraliser. Note that the number of conjugacy 
classes of maximal tori in $G$ is bounded by the number of different cyclotomic 
polynomials evaluated in $q$ appearing as factors of $|G|$. Moreover the power 
of a cyclotomic polynomial in the order of $G$ gives the rank of the 
corresponding maximal torus (except possibly when the prime divides the order 
of the centre but in this case the Sylow subgroup is not cyclic, 
see~\cite[(33.14)]{A})

Recall $\Sigma$ is the root system associated to $G$ as in \cite[2.3.1]{GLS}; 
let $\Pi=\{\alpha_1, \ldots , \alpha_n\}$ be a fundamental system for
$\Sigma$ as in Table~1.8 in~\cite{GLS}, $\alpha_{\ast}$ be the lowest root
relative to $\Pi$ as defined in~\cite[Paragraph~1.8]{GLS} and set
$\Pi^{\ast}=\Pi \cup \{\alpha_{\ast}\}$. We recall that $|G|$ can be deduced 
from \cite[Table 2.2]{GLS} and the Dynkin diagrams can be found in
\cite[Table 1.8]{GLS}. Observe that, by \cite[Proposition~2.6.2]{GLS}, if 
$\Sigma_0$ is a root subsystem of $\Sigma$, then $G$ contains a subsystem 
subgroup $H$, which is a central product of groups of Lie type corresponding to 
the irreducible constituents of $\Sigma_0$. In order to prove the lemma, we 
shall show that for every group $G$ and for every element $g$ of order a prime 
$r$ lying in a maximal torus belonging to any but four conjugacy classes of 
maximal tori, either the Sylow $r$-subgroup is not cyclic or we find a 
subsystem subgroup $H$ that is a central product of two groups $H_1$ and $H_2$ 
such that $H_1$ contains $g$ and $H_2$ is not abelian. Note that for every 
prime power $q$, $A_1(q)$ is a non-abelian group (see
\cite[Theorem~2.2.7]{GLS}).

\medskip
We treat the case $G\cong A_n(q)$ in detail as an example. All other cases can 
be dealt similarly. Assume $G\cong A_n(q)$. Let $m$ be the minimum index $i$ 
such that $r$ divides $q^{i+1}-1$ and let $\Sigma_0$ be generated by 
$\Pi^{\ast}\setminus \{\alpha_1,\alpha_{n}\}$. Then the corresponding subsystem 
subgroup is $H=H_1 \cdot H_2$, where $H_1\cong A_{n-2}(q)$ and 
$H_2\cong A_1(q)$. Thus if $m\leq n-1$, then $H_1$ contains an element $g$ of 
order $r$ and $C_G(g)$ contains $H_2$ which is not abelian. Therefore, since 
$g$ has an abelian centraliser, $r$ may divide only $(q^n-1)(q^{n+1}-1)$, that 
is $r$ divides $\Phi_{n}(q)\Phi_{n+1}(q)$. Hence we have at most two conjugacy 
classes of maximal tori with abelian centraliser.
\qed

%%%%%%%%%%%%%%%%%%%%%%%%%%%%%%%%%%%%%%%%%%%%%%%%%%%%%%%%%%%%

\section{Proof of Theorem~\ref{th:main4} }\label{sec:proof}

Let $G$ be a finite group of orientation preserving diffeomorphisms of $M$, a 
closed orientable connected $3$-manifold which is not homeomorphicto $\S^3$. 
%We will prove that $G$ contains at most six conjugacy classes of cyclic 
%subgroups generated by a hyperelliptic rotation of order not a power of $2$. 

As noted in Section~\ref{sec:road-map}, there at most nine conjugacy classes of 
cyclic subgroups generated by a hyperelliptic rotation of order $2^n$, so we 
concentrate on hyperelliptic subgroups whose order is not a power of two and we 
will prove that there are at most six conjugacy classes of such subgroups.

Let $\mathcal{S}=\{C_1,\dots,C_m\}$ be the set of cyclic subgroups of odd prime 
orders that are generated by powers of the hyperelliptic rotations of $G$. We 
recall that the conjugate of a power of a hyperelliptic rotation is the power 
of a hyperelliptic rotation too. By Remark~\ref{rem:normaliser} and 
Corollary~\ref{cor:2hyprot}, $\mathcal{S}$ is an algebraically hyperelliptic 
collection, and actually the definition of an algebraically hyperelliptic 
collection was chosen precisely to capture the behaviour of cyclic subgroups of odd prime orders generated by powers of the hyperelliptic rotations. For 
$i\in \{1, \ldots , m\}$, let $p_i$ be the prime order of $C_i$ and $\pi$ be 
the set including every $p_i$. We denote by $G_0$ the subgroup generated by the 
subgroups $C_i.$ Let $A$ be the maximal normal solvable subgroup of $G_0$ of 
order coprime with all $p_i$. We denote by $\overline{G_0}$ the quotient group 
$G_0/A$ and by $\overline C_i$ the projection of $C_i$ to $\overline{G_0}$.
\smallskip

\textbf{Case 1.} If $G_0$ is solvable, we are in the hypotheses of 
Lemma~\ref{lem:algebraic-solvable} which affirms that up to conjugacy the 
hyperelliptic rotations of $G$ commute. Since they commute we can apply 
Lemma~\ref{lem:abelian-case} and conclude that there are at most three 
conjugacy classes of subgroups generated by hyperelliptic rotations of order 
not a power of two. 

Even if its proof is easier than that of the remaining situations, the 
solvable case is interesting in its own right, in particular for it plays an 
important role also in the proof of Theorem~\ref{th:knots}. We then summarise 
the conclusions in the following proposition.

\begin{Proposition}\label{pro:solvable}
Let $M$ be a closed $3$-manifold not homeomorphic to $\S^3$. Let $G$ be a 
finite group of orientation preserving diffeomorphisms of $M$. If $G$ is 
solvable, it contains at most three conjugacy classes of hyperelliptic 
subgroups of order not a power of two, and any two such subgroups commute up to 
conjugacy. Moreover, if there are three conjugacy classes, then their orders
must be pairwise coprime.
\end{Proposition}

\smallskip

\textbf{Case 2.} 
Suppose $G_0$ is not solvable and it has no rotation of order $2$ outside $A$. 
Then, because of the structure of the normaliser of a hyperelliptic rotation as 
described in Remark~\ref{rem:normaliser}, $G_0$ has no involution acting 
dihedrally on any $C_i$ and by Part 3 of 
Theorem~\ref{lem:algebraic-non-solvable} describing the structure of $G_0$, 
$E(\overline{G_0})$ is simple and, for each $p_i$, any Sylow $p_i$-subgroup 
contains exactly two distinct conjugates of $\bar C_i$. By 
Remark~\ref{rem:normaliser} every hyperelliptic rotation commuting with one of 
these two subgroups of order $p_i$ commutes also with the other one. From this 
fact, Lemma~\ref{lem:algebraic-solvable}, and Lemma~\ref{lem:abelian-case}, it 
follows that $\gamma^s_{\pi}(G_0)$ bounds from above the number of conjugacy 
classes of hyperelliptic subgroups of order not a power of two (see also 
Proposition~\ref{pro:solvable}). It is easy to see that 
$\gamma^s_{\pi}(G_0)\leq \gamma^s_{\pi}(E(\overline{G_0})).$ By 
Proposition~\ref{simple} we have $\gamma^s_{\pi}(E(\overline{G_0}))\le 4$ and so
we get the thesis in this case.

\smallskip

\textbf{Case 3.} Suppose $G_0$ is not solvable and it has a rotation of order 
$2$ not contained in $A$. The groups containing a rotation of order two are 
studied in \cite{Mec} where the following result was proved. 

\begin{Theorem}\label{thm:involution}\cite{Mec}
Let $D$ be a finite group of orientation-preserving diffeomorphisms of a closed 
orientable $3$-manifold. Let $O$ be the maximal normal subgroup of odd order 
and $E(\tilde{D})$ be the layer of $\tilde{D}=D/O$. Suppose that $D$ contains 
an involution which is a rotation.

\begin{enumerate}

\item If $E(\tilde {D})$ is trivial, there exists a normal subgroup $H$ of $D$ 
such that $H$ is solvable and $D/H$ is isomorphic to a subgroup of 
$\mathbb{A}_8$, the alternating group on $8$ letters.

\item If the semisimple group $E(\tilde D)$ is not trivial, it has at most two 
components and the factor group of $\tilde D/E(\tilde  D)$ is solvable. 
%Moreover, the factor group of $E(\tilde {G})$ by its centre is either a simple 
%group of sectional $2$-rank at most $4$ or the direct product of two simple 
%groups of sectional $2$-rank at most $2$.

\end{enumerate}

\noindent
Moreover if $D$ contains a rotation of order $2$ such that its projection is 
contained in $E(\tilde D)$, then $E(\tilde D)$ is isomorphic to one of the 
following groups:
$$PSL_2(q), \;\;   SL_2(q)\times_{\ZZ/2} SL_2(q')$$
where $q$ and $q'$ are odd prime powers greater than $4$. 
\end{Theorem}

Applying Theorem~\ref{thm:involution} and 
Theorem~\ref{lem:algebraic-non-solvable} to $G_0$, we get that 
$E(\overline{G_0})$ is isomorphic either to a subgroup of $\mathbb{A}_8$, or to 
$PSL_2(q)$, or to $PSL_2(q)\times PSL_2(q')$. In the first case, there are at 
most three odd primes dividing the order of $E(\overline{G_0})$ and the thesis 
follows again from Lemma~\ref{lem:abelian-case} and 
Lemma~\ref{lem:algebraic-solvable}. 

In the remaining cases, we will use a solvable normal $\pi$-cover to bound the 
number of conjugacy classes. We have that $\gamma^s_{\pi}(PSL_2(q))\leq 2.$ In 
fact the upper triangular matrices form a solvable subgroup of $SL_2(q)$ of 
order $(q-1)q$, moreover $SL_2(q)$ contains a cyclic subgroup of order $q+1$ 
(see \cite{H}). The conjugates of the projections of these two subgroups to 
$PSL_2(q)$ give a solvable normal $\pi$-cover of $PSL_2(q)$. It is again easy 
to see that, if $E(\overline{G_0})$ is isomorphic to $PSL_2(q)$, then 
$\gamma^s_{\pi}(G)\leq 2$. As above, by Lemma~\ref{lem:algebraic-solvable} and  
Lemma~\ref{lem:abelian-case}, we get a bound of six in this case. 

Finally, if $E(\overline{G})$ is isomorphic to $PSL_2(q)\times PSL_2(q')$, then 
it follows from the discussion above that 
$\gamma^s_{\pi}(E(\overline{G_0}))\leq 4$, and hence 
$\gamma^s_{\pi}(G_0)\leq 4.$ Now, by Part 1 of 
Theorem~\ref{lem:algebraic-non-solvable}, for each $p_i\in \pi$, $p_i$ divides
the order of both $PSL_2(q)$ and $PSL_2(q')$, so that the Sylow $p_i$-subgroup 
of $\overline{G_0}$ has rank $2$. For each Sylow $p_i$-subgroup $S_i$, one can 
find in $E(\overline{G_0})$ elements of the group which normalise $S_i$ but do 
not act in the same way on all of its elements. Indeed, let $H$ be a cyclic 
subgroup of order $p_i$ in $PSL_2(q)$. If $p_i$ does not divide $q$, then there 
is an element of order $2$ in $PSL_2(q)$ which acts dihedrally on $H$ but 
commutes with all elements of order $p$ in $PSL_2(q')$. If $p_i$ divides $q$, 
then $q=p_i$ since the Sylow $p_i$-subgroup of the component must be cyclic. 
Since $q>3$, we have again elements in $PSL_2(q)$ that normalise the Sylow 
$p$-subgroup but do not centralise it (note that in this case the structure of 
the normaliser is compatible with the description given in
Remark~\ref{rem:normaliser} only if $p=q=5$). We now deduce from 
Remark~\ref{rem:conjugate-rotations} and Corollary~\ref{cor:powers suffice}
that $S_i$ contains two subgroups conjugate to $C_i$.

%As in Case 2, using Remark~\ref{rem:conjugate-rotations}, one can see that any 
%Sylow $p_i$-subgroup contains exactly two distinct and conjugate cyclic 
%subgroups which are generated by the power of a hyperelliptic rotation and thus 
Reasoning as in Case 2 we get that $\gamma^s_{\pi}(G_0)$ bounds from above the 
number of conjugacy classes of hyperelliptic subgroups of order not a power of 
two. This concludes the proof. 
\qed

\begin{remark}\label{r:casesbysymmetry} 
We wish to stress that the three cases that appear in the proof of
Theorem~\ref{th:main4} correspond to different types of symmetries possessed by
the knots that are branched covered by the manifold. Let $M$ be a
closed, connected, oriented $3$-manifold and let $G$ a finite group of
diffeomorphisms of $M$ generated by a set $\{\psi_1, \dots, \psi_k\}$ of 
hyperelliptic rotations. Denote by $K_i$ the quotient knot
$Fix(\psi_i)/\psi_i$, $i=1,\dots, k$. If there exists a $K_i$ which is neither
strongly invertible nor self-symmetric with respect to $G$, then we are in Case
1, that is, the hyperelliptic rotations commute up to conugacy (see \cite[Thm
2.10, page 144]{Su}). If there is a $K_i$ which is strongly invertible with 
respect to $G$, then $G$ contains a rotation of order $2$ and we are in Case 3. 
Otherwise, every $K_i$ is self-symmetric with respect to $G$ but not strongly 
invertible with respect to $G$, and we are in Case 2.
\end{remark}

The proof of Theorem~\ref{th:main4} shows that topology and geometry can impose
extra constraints on the conditions that can be derived by the algebra alone.
In this spirit, note that if $M$ is a (closed, connected, oriented) reducible 
$3$-manifold and $G$ a finite group of diffeomorphisms of $M$ generated by
hyperelliptic rotations, then, by the equivariant sphere theorem, $G$ is
isomorphic to a finite subgroup of $SO(3)$, that is cyclic, dihedral or a
spherical triangular group. It follows readily, that, up to conjugacy, at most
three cyclic hyperelliptic subgroups can be contained in any such $G$. We know,
however, that the group of diffeomorphisms of a reducible $3$-manifolds can 
admit arbitrarily many conjugacy classes of cyclic hyperelliptic subgroups. 

The main point here is that, generically, one expects that two hyperelliptic 
rotations in the group $Diff^+(M)$ generate an infinite subgroup.
%a set of hyperelliptic rotations will
%not generate a finite group of diffeomorphisms. 
As a consequence, 
Theorem~\ref{th:main4} cannot be directly exploited to obtain bounds for 
non-hyperbolic manifolds and new strategies must be developped, as we will see
in the next section.

\section{Proof of Theorem \ref{th:knots}}\label{sec:knots}

The statement of Theorem \ref{th:knots} is equivalent to the following:

\begin{Theorem}\label{th:main2}
Let $M$ be a closed, orientable, connected, irreducible $3$-manifold which is
not homeomorphic to $\S^3$, then the group $Diff^+(M)$ of orientation
preserving diffeomorphisms of $M$ contains at most six conjugacy classes of
hyperelliptic subgroups of odd prime order.
\end{Theorem}

%Notice that generically one expects that two hyperelliptic rotations in the 
%group $Diff^+(M)$ generate an infinite subgroup.

\subsection{Proof of Theorem \ref{th:main2} for Seifert manifolds}
\label{sec:seifert} 

In this section we prove Proposition~\ref{prop:seifert} which implies 
Theorem~\ref{th:main2} for closed Seifert fibred $3$-manifolds. We also show 
that the assumption that the hyperelliptic rotations have odd prime orders 
cannot be avoided in general by exhibiting examples of closed Seifert fibred 
$3$-manifolds $M$ such that $Diff^+(M)$ contains an arbitrarily large number of 
conjugacy classes of hyperelliptic subgroups of 
odd, but not prime, orders.

\begin{Proposition}\label{prop:seifert}
Let $M$ be a closed Seifert fibred $3$-manifold which is not homeomorphic to 
$\S^3$. Then the group $Diff^+(M)$ of orientation preserving diffeomorphisms of 
$M$ contains at most one conjugacy class of hyperelliptic subgroups of odd 
prime order except if $M$ is a Brieskorn integral homology sphere with $3$ 
exceptional fibres. In this latter case $Diff^+(M)$ contains at most three non 
conjugate hyperelliptic subgroups of odd prime orders.
\end{Proposition}

\demo

By hypothesis $M$ is a cyclic cover of $\S^3$ branched over a knot, so 
it is orientable and a rational homology sphere by Remark~\ref{rem:Qsphres}. 
Notably, $M$ cannot be $\S^1 \times \S^2$ nor a Euclidean manifold, except 
for the Hantzsche-Wendt manifold, see \cite[Chap. 8.2]{Or}. In particular, 
since $M$ is prime it is also irreducible.

Consider a hyperelliptic rotation $\psi$ on $M$ of odd prime order $p$ and let 
$K$ be the image of $Fix(\psi)$ in the quotient $\S^3= M/\psi$ by the action of 
$\psi$. The knot $K$ must be hyperbolic or a torus knot, otherwise its exterior 
would be toroidal and have a non-trivial JSJ-collection of essential tori 
which would lift to a non-trivial JSJ-collection of tori for $M$, since the 
order of $\psi$ is $p >2$ (see \cite{JS,J} and \cite{BS}). By the orbifold 
theorem (see \cite{BoP}, \cite{CHK}), the cyclic branched cover with order 
$p \geq 3$ of a hyperbolic knot is hyperbolic, with a single exception for 
$p = 3$ when $K$ is the the figure-eight knot and $M$ is the Hantzsche-Wendt 
Euclidean manifold. But then, by the orbifold theorem and the classification of 
$3$-dimensional christallographic groups, $\psi$ generates the unique, up to 
conjugacy, Euclidean hyperelliptic cyclic subgroup of $Diff^+(M)$, see 
for example \cite{Dun}, \cite{Z1}.

So we can assume that $M$ is the $p$-fold cyclic cover of $\S^3$ branched along 
a non-trivial torus knot $K$ of type $(a,b)$, where $a > 1$ and $b > 1$ are 
coprime integers. Then $M$ is a Brieskorn-Pham manifold 
$M= V(p, a, b) = \{z^p + x^a + y^b = 0 \, \, \, \text{with}  \, \, \,  
(z, x, y) \in \CC^3 \, \, \,  \text{and} \, \, \, \vert z \vert^2 + 
\vert x \vert ^2 + \vert y \vert^2 = 1 \}$. 
A simple computation shows that $M$ admits a Seifert fibration with $3$, $p$ or 
$p+1$ exceptional fibres and base space $\S^2$, see \cite[Lem. 2]{Ko}, or 
\cite[Lemma 6 and proof of Lemma 7]{BoPa}. In particular $M$ has a unique 
Seifert fibration, up to isotopy: by \cite{Wa}, \cite{Sco} and \cite{BOt} the 
only possible exception with base $\S^2$ and at least $3$ exceptional fibres 
is the double of a twisted $I$-bundle, which is not a rational homology sphere, 
since it fibers over the circle. We distinguish now two cases:\\

\noindent{\bf Case 1}: \emph{The integers $a$ and $b$ are coprime with $p$, and 
there are three singular fibres of pairwise relatively prime orders $a$, $b$ 
and $p$}. By the orbifold theorem any hyperelliptic rotation of $M$ of order 
$>2$ is conjugate into the circle action $\S^1 \subset Diff^{+}(M)$ inducing the 
Seifert fibration, hence the uniqueness of the Seifert fibration, up to 
isotopy, implies that $M$ admits at most $3$ non conjugate hyperelliptic groups 
of odd prime orders belonging to 
$\{ a, b, p\}$. Indeed $M$ is a Brieskorn integral homology sphere, see 
\cite{BPZ}.\\

\noindent{\bf Case 2}: \emph{Either $a=p$ and $M$ has $p$ singular fibres of 
order $b$, or $a = a'p$ with $a' > 1$, and $M$ has $p$ singular fibres of order 
$b$ and one extra singular fibre of order $a'$.} In both situations, there are 
$p \geq 3$ exceptional fibres of order $b$ which are cyclically permuted by the 
hyperelliptic rotation $\psi$. As before, $M$ has a unique Seifert fibration, 
up to isotopy. Therefore, up to conjugacy, $\psi$ is the only hyperelliptic 
rotation of order $p$ on $M$, and by the discussion above $M$ cannot admit a 
hyperelliptic rotation of odd prime order $q \neq p$. 
\qed

\begin{remark}\label{rem:brieskorn}
The requirement that the rotations are hyperelliptic is essential in the proof 
of Proposition~\ref{prop:seifert}. The Brieskorn homology sphere  
$\Sigma(p_1,\ldots,p_n), \, \, n \geq 4$, with $n \geq 4$ exceptional fibres 
admits $n$ rotations of pairwise distinct prime orders but which are not 
hyperelliptic.

The hypothesis that the orders of the hyperelliptic rotations are $\neq 2$ 
cannot be avoided either.

Indeed, Montesinos's construction of fibre preserving hyperelliptic involutions 
on Seifert fibered rational homology spheres \cite{Mon1}, \cite{Mon2}, (see 
also \cite[Appendix A]{BS}, \cite[Chapter 12]{BZH}), shows that for any given 
integer $n$ there are infinitely many closed orientable Seifert fibred 
$3$-manifolds with at least $n$ conjugacy classes of hyperelliptic rotations of 
order $2$.

On the other hand, the hypothesis that the orders are odd primes is sufficient
but not necessary: A careful analysis of the Seifert invariants shows that if
$M\neq\S^3$ is a Seifert rational homology sphere, then $M$ can be the cyclic
branched cover of a knot in $\S^3$ of order $>2$ in at most three ways.

The hypotheses of Proposition~\ref{prop:seifert} cannot be relaxed further, 
though: Proposition~\ref{prop:circle bundles} below shows that there exist 
closed $3$-dimensional circle bundles with arbitrarily many conjugacy classes 
of hyperelliptic rotations of odd, but not prime, orders. 
\end{remark}

\begin{Proposition}\label{prop:circle bundles} 
Let $N$ be an odd prime integer. For any integer $1 \leq q < \frac{N}{2}$ the 
Brieskorn-Pham manifold 
$M = V((2^q + 1)(2^{(N-q)} + 1), 2^q + 1, 2^{(N-q)} + 1)$ is a circle bundle 
over a closed surface of genus $g = 2^{N - 1}$ with Euler class $\pm1$. Hence, 
up to homeomorphism (possibly reversing the orientation), $M$ depends only on 
the integer $N$ and admits at least $\frac{N-1}{2}$ conjugacy classes of hyperelliptic groups of odd orders.
\end{Proposition}

\demo
We remark first that the integers $q$ and $N-q$ are relatively prime, because 
$N$ is prime. If $k$ is a common prime divisor of $2^q + 1$ and 
$2^{(N-q)} + 1$, by the Bezout identity we have 
$2^1=2^{aq+b(N-q)}\equiv (-1)^{a+b}\mod k$ which implies that $k=3$. But then 
$(-1)^q \equiv (-1)^{(N-q)} \equiv -1 \mod 3$ and thus $(-1)^N \equiv 1 \mod 3$ 
which is impossible since $N$ is odd. Hence $2^q+ 1$ and $2^{(N-q)} + 1$ are 
relatively prime. 

So the Brieskorn-Pham manifold $M$ is the $(2^q+ 1)(2^{(N-q)} + 1)$-fold cyclic 
cover of $\S^3$ branched over the torus knot $K_q = T(2^q + 1, 2^{(N-q)} + 1)$. 
It is obtained by Dehn filling the $(2^q + 1)(2^{(N-q)} + 1)$-fold cyclic cover 
of the exterior of the torus knot $K_q$ along the lift of its meridian. The 
$(2^q + 1)(2^{(N-q)} + 1)$-fold cyclic cover of the exterior of $K_q$ is a 
trivial circle bundle over a once punctured surface of genus $g = 2^N - 1$. On 
the boundary of the torus-knot exterior the algebraic intersection between
a meridian and a fibre of the Seifert fibration of the exterior is $\pm 1$ 
(the sign depends on a choice of orientation, see for example 
\cite[Chapter 3]{BZH}). So on the torus boundary of the 
$(2^q + 1)(2^{(N-q)} + 1)$-fold cyclic cover the algebraic intersection between 
the lift of a meridian of the torus knot and an $\S^1$-fiber is again $\pm 1$.  
Hence the circle bundle structure of the $(2^q + 1)(2^{(N-q)} + 1)$-fold cyclic 
cover of the exterior of the torus knot $K_q$ can be extended with Euler class 
$\pm 1$ to the Dehn filling along the lift of the meridian. So $M$ is a circle 
bundle over a closed surface of genus $g = 2^{N - 1}$ with Euler class $\pm 1$. 

Since the torus knots $K_q = T(2^q + 1, 2^{(N-q)} + 1)$ are pairwise 
inequivalent for $1 \leq q \leq\frac{N-1}{2}$, the hyperelliptic subgroups 
corresponding to the $(2^q+ 1)(2^{(N-q)} + 1)$-fold cyclic branched covers of 
the knots $K_q$ are pairwise not conjugate in $Diff^{+}(M)$. 
\qed

\begin{remark}\label{r:not finite}
Note that the Seifert manifolds $M$ and their hyperelliptic rotations 
constructed in Proposition~\ref{prop:circle bundles} enjoy the following 
properties: If $N>8$, then no hyperelliptic rotation can commute up to 
conjugacy with all the remaining ones (see Proposition~\ref{pro:solvable} and 
\cite[Theorem 2]{BoPa}). If $N>14$ no finite subgroup of
$Diff^+(M)$ can contain up to conjugacy all hyperelliptic rotations of $M$, 
according to Theorem~\ref{th:main4}.
\end{remark}

%%%%%%%%%%%%%%%%%%%%%%%%%%%%%%%%%%%%%%%%%%%%%%%%%%%%%%%%%%%%%%%%%%%%%%%%%%%%%%%%

\subsection{Reduction to the finite group action case}\label{sec:reduction}

The fact that Theorem~\ref{th:main2} implies Corollary~\ref{th:main} follows 
from the existence of a decomposition of a closed, orientable $3$-manifold $M$ 
as a connected sum of prime manifolds and the observation that a hyperelliptic 
rotation on $M$ induces a hyperelliptic rotation on each of its prime summands. 
A $3$-manifold admitting a hypereliptic rotation must be a rational homology 
sphere, and so $M$ cannot have $\S^2\times\S^1$ summands. Hence all prime 
summands are irreducible and at least one is not homeomorphic to $\S^3$, since 
$M$ itself is not homeomorphic to $\S^3$. This is enough to conclude.\\

The remaining of this section is devoted to the proof that 
Theorem~\ref{th:main4} implies Theorem~\ref{th:main2}.

We prove the following proposition:

\begin{Proposition}\label{prop:finitereduction}
If $M$ is a closed, orientable, irreducible $3$-manifold such that there are 
$k \geq 7$ conjugacy classes of hyperelliptic subgroups of $Diff^{+}(M)$ whose 
order is an odd prime, then $M$ is homeomorphic to $\S^3$.
\end{Proposition}

\demo

Let $M$ be a closed, orientable, irreducible $3$-manifold such that
$Diff^{+}(M)$ contains $k \geq 7$ conjugacy classes of hyperelliptic subgroups of  odd prime orders.

According to the orbifold theorem (see \cite{BoP}, \cite{BMP}, \cite{CHK}), a
closed orientable irreducible manifold $M$ admitting a rotation has geometric
decomposition. This means that $M$ can be split along a (possibly empty) finite
collection of $\pi_1$-injective embedded tori into submanifolds carrying either
a hyperbolic or a Seifert fibered structure. This splitting along tori is
unique up to isotopy and is called the JSJ-decomposition of $M$, see for 
example \cite{NS}, \cite[chapter 3]{BMP}. In particular, if its 
JSJ-decomposition is trivial, $M$ admits either a hyperbolic or a Seifert 
fibred structure.

First we see that $M$ cannot be hyperbolic. Indeed, if the manifold $M$ is 
hyperbolic then, by the orbifold theorem, any hyperelliptic rotation is 
conjugate into the finite group $Isom^+(M)$ of orientation preserving 
isometries of $M$. Hence, applying Theorem~\ref{th:main4} to $G = Isom^+(M)$, 
we see that $k\le 6$ against the hypothesis.
 
If the manifold $M$ is Seifert fibred, it follows readily from 
Proposition~\ref{prop:seifert} of the previous section that $M=\S^3$. So we are 
left to exclude the case where the JSJ-decomposition of $M$ is not empty.

Consider the JSJ-decomposition of $M$: each geometric piece admits either a 
complete hyperbolic structure with finite volume or a Seifert fibred product 
structure with orientable base. Moreover, the geometry of each piece is 
unique, up to isotopy.

Let $\Psi = \{\psi_1,\dots, \psi_k, \, k \geq 7\}$ be the set of hyperelliptic 
rotations which generate non conjugate cyclic subgroups in $Diff^+(M)$. By the 
orbifold theorem \cite{BoP}, \cite{BMP}, \cite{CHK}, after conjugacy, one can 
assume that each hyperelliptic rotation preserves the JSJ-decomposition,
acts isometrically on the hyperbolic pieces, and respects the product structure 
on the Seifert pieces. We say that they are \emph{geometric}.

Let $\Gamma$ be the dual graph of the JSJ-decomposition: it is a tree, for 
$M$ is a rational homology sphere (in fact, the dual graph of the
JSJ-decomposition for a manifold which is the cyclic branched cover of a knot 
is always a tree, regardless of the order of the covering). Let 
$H \subset Diff^+(M)$ be the group of diffeomorphisms of $M$ generated by the 
set $\Psi$ of geometric hyperelliptic rotations. By \cite[Thm 1]{BoPa}, there 
is a subset $\Psi_0 \subset \Psi$ of $k_0 \geq 4$ hyperelliptic rotations with 
pairwise distinct odd prime orders, say 
$\Psi_0 =\{\psi_i, \, i = 1,\dots,k_0 \}$.

Let $H_\Gamma$ denote the image of the induced representation of $H$ in 
$Aut(\Gamma)$. Since rotations of finite odd order cannot induce an inversion 
on any edge of $\Gamma$, the finite group $H_\Gamma$ must fix pointwise a 
non-empty subtree $\Gamma_{f}$ of $\Gamma$.  

The idea of the proof is now analogous to the ones in \cite{BoPa} and 
\cite{BPZ}: we start by showing that, up to conjugacy, the $k_0 \geq 4$ 
hyperelliptic rotations with pairwise distinct odd prime orders can be chosen 
to commute on the submanifold $M_f$ of $M$ corresponding to the subtree 
$\Gamma_f$. We consider then the maximal subtree corresponding to a submanifold 
of $M$ on which these hyperelliptic rotations commute up to conjugacy and 
prove that such subtree is in fact $\Gamma$. Then the conclusion follows 
%as in the proof of Theorem~\ref{th:main4}, 
by applying Lemma~\ref{lem:abelian-case}.
%\ref{prop:commuting rotations}.

\medskip 

The first step of the proof is achieved by the following proposition:

\begin{Proposition}\label{prop:commuting submanifold}
The hyperelliptic rotations in $\Psi_0$ commute, up to conjugacy in 
$Diff^+(M)$, on the submanifold $M_f$  of $M$ corresponding to the subtree 
$\Gamma_f$.
\end{Proposition}

\demo

Since the hyperelliptic rotations in $\Phi$ have odd orders, either $\Gamma_f$
contains an edge, or it consists of a single vertex. We shall analyse these two
cases. 

\medskip

\noindent{\bf Case 1: $M_f$ contains an edge.}

\begin{Claim}\label{claim:invariant edge}
Assume that $\Gamma_f$ contains an edge and let $T$ denote the corresponding
torus. The hyperelliptic rotations in $\Psi_0$ commute, up to conjugacy in 
$Diff^+(M)$, on the geometric pieces of $M$ adjacent to $T$.
\end{Claim}

\demo

The geometric pieces adjacent to $T$ are left invariant by the hyperelliptic 
rotations in $\Psi_0$, since their orders are odd. Let $V$ denote one of the 
two adjacent geometric pieces: each hyperelliptic rotation acts non-trivially 
on $V$ with odd prime order. We distinguish two cases according to the geometry 
of $V$.

\noindent {\bf $V$ is hyperbolic.} In this case all rotations act as isometries 
and leave a cusp invariant. Since their order is odd, the rotations must act as 
translations along horospheres, and thus commute. 

Note that, even in the case where a rotation has order $3$, its axis cannot 
meet a torus of the JSJ-decomposition of $M$ for each such torus is separating 
and cannot meet the axis in an odd number of points.

\noindent {\bf $V$ is Seifert fibred.} In this case we can assume that the 
hyperelliptic rotations in $\Psi$ preserve the Seifert fibration with 
orientable base. Since their orders are odd and prime, each one preserves the 
orientation of the fibres and of the base. The conjugacy class of a 
fiber-preserving rotation of $V$ with odd prime order depends only on its 
combinatorial behaviour, i.e. its translation action along the fibre and the 
induced permutation on cone points and boundary components of the base. In 
particular, two geometric rotations with odd prime order having the same 
combinatorial data are conjugate via a diffeomorphism isotopic to the identity.

Since the hyperelliptic rotations in $\Psi_0$ have pairwise distinct odd prime 
orders, an analysis of the different cases described in Lemma~\ref{lem:seifert} 
below shows that at most one among these hyperelliptic rotations can induce a 
non-trivial action on the base of the fibration, and thus the remaining ones 
act by translation along the fibres and induce the identity on the base. Since 
the translation along the fibres commutes with every fiber-preserving 
diffeomorphism of $V$, the hyperelliptic rotations in $\Psi_0$ commute on $V$.

Lemma~\ref{lem:seifert} describes the Seifert fibred pieces of a manifold 
admitting a hyperelliptic rotation of odd prime order, as well as the action of 
the rotation on the pieces. Its proof can be found in \cite[Lemma 6 and proof 
of Lemma 7]{BoPa}, see also \cite[Lem. 2]{Ko}.

\begin{Lemma}\label{lem:seifert}
Let $M$ be an irreducible $3$-manifold admitting a non-trivial
JSJ-decomposition. Assume that $M$ admits a hyperelliptic
rotation of prime odd order $p$. Let $V$ be a Seifert piece of the 
JSJ-decomposition for $M$. Then the base $B$ of $V$ can be:

\begin{enumerate}

\item A disc with $2$ cone points. In this case either the rotation freely 
permutes $p$ copies of $V$ or leaves $V$ invariant and acts by translating 
along the fibres.

\item A disc with $p$ cone points. In this case the rotation leaves $V$ 
invariant and cyclically permutes the singular fibres, while leaving a regular
one invariant. 

\item A disc with $p+1$ cone points. In this case the rotation leaves $V$ and a 
singular fibre invariant, while cyclically permuting the remaining $p$ singular
fibres.

\item An annulus with $1$ cone point. In this case either the rotation freely 
permutes $p$ copies of $V$ or leaves $V$ invariant and acts by translating 
along the fibres.

\item An annulus with $p$ cone points. In this case the rotation leaves $V$
invariant and cyclically permutes the $p$ singular fibres.

\item A disc with $p-1$ holes and $1$ cone point. In this case the rotation
leaves $V$ invariant and cyclically permutes all $p$ boundary components, while
leaving invariant the only singular fibre and a regular one.

\item A disc with $k$ holes, $k\ge2$. In this case either the rotation freely
permutes $p$ copies of $V$ or leaves $V$ invariant. In this latter case either 
the rotation acts by translating along the fibres, or $k=p-1$ and the rotation
permutes all the boundary components (while leaving invariant two fibres), or 
$k=p$ and the rotation permutes $p$ boundary components, while leaving 
invariant the remaining one and a regular fibre. 
\end{enumerate} 
\end{Lemma}
\qed
\smallskip

We conclude that the rotations in $\Psi_0$ can be chosen to commute on the 
submanifold $M_f$ of $M$ corresponding to $\Gamma_f$ by using inductively at 
each edge of $\Gamma_f$ the gluing lemma below (see [Lemma 6]\cite{BPZ}). We 
give the proof for the sake of completeness.

\begin{Lemma}\label{lem:gluing} 
If the rotations preserve a JSJ-torus $T$ then they commute on the union of 
the two geometric pieces adjacent to $T$.
\end{Lemma}

\demo

Let $V$ and $W$ be the two geometric pieces adjacent to $T$. By Claim 
\ref{claim:invariant edge}, after conjugacy in $Diff^{+}(M)$, the rotations in 
$\Psi_0$ commute on $V$ and $W$. Since they have pairwise distinct odd prime 
orders, their restrictions on $V$ and $W$ generate two cyclic groups of the 
same finite odd order. Let $g_V$ and $g_W$ be generators of these two cyclic 
groups. Since they have odd order, they both act by translation on $T$. We need 
the following result about the slope of translation for such periodic 
transformation of the torus:

\begin{Claim}\label{claim:slope} 
Let $\psi$ be a periodic diffeomorphism of the product $T^2 \times [0,1]$ which 
is isotopic to the identity and whose restriction to each boundary torus 
$T\times \{i\}$, $i= 0, 1,$ is a translation with rational slopes $\alpha_0$ 
and $\alpha_1$ in $H_1(T^2; \ZZ)$. Then $\alpha_0 = \alpha_1$.
\end{Claim}

\demo 

By Meeks and Scott \cite[Thm 8.1]{MS}, see also \cite[Prop. 12]{BS}, there is a 
Euclidean product structure on $T^2 \times [0,1]$ preserved by $\psi$ such that 
$\psi$ acts by translation on each fiber $T \times\{t\}$ with rational slope 
$\alpha_t$. By continuity, the rational slopes $\alpha_t$ are constant.
\qed

\medskip

Now the the following claim shows that the actions of $g_W$ and $g_V$ can be 
glued on $T$.

\begin{Claim}\label{claim:identification}
The translations $g_{V}\vert_T$ and $g_{W}\vert_T$ have the same slope in 
$H_1(T^2; \ZZ)$. 
\end{Claim}

\demo 

Let $\Psi_0 =\{\psi_i, \, i = 1,\dots,k_0 \}$. Let $p_i$ the order of $\psi_i$ 
and $q_i = \Pi_{j \neq i} p_j$. Then the slopes $\alpha_V$ and $\alpha_W$ of 
$g_{V}\vert_T$ and $g_{W}\vert_T$ verify: $q_i \alpha_V = q_i \alpha_W$ for 
$i= 1,...,k_0$, by applying Claim \ref{claim:slope} to each $\psi_i$. Since the 
$GCD$ of the $q_i$ is $1$, it follows that $\alpha_V = \alpha_W$.
\qed

\medskip

This finishes the proof of Lemma~\ref{lem:gluing} and of 
Proposition~\ref{prop:commuting submanifold} when $M_f$ contains an edge. 
\qed

\medskip

To complete the proof of Proposition~\ref{prop:commuting submanifold} it 
remains to consider the case where $\Gamma_f$ is a single vertex.

\medskip

\noindent{\bf Case 2: $M_f$ is a vertex.}

\medskip

\begin{Claim}\label{claim:one vertex}
Assume that $\Gamma_f$ consists of a single vertex and let $V$ denote the 
corresponding geometric piece. Then the hyperelliptic rotations in $\Psi_0$ 
commute on $V$, up to conjugacy in $Diff^+(M)$.
\end{Claim}

\demo

We consider again two cases according to the geometry of $V$. 

The case where $V$ is Seifert fibred follows once more from 
Lemma~\ref{lem:seifert}. 

We consider now the case where $V$ is hyperbolic.

In this case, the hyperelliptic rotations in $\Psi$ act non-trivially on $V$ 
by isometries of odd prime orders. The restriction $H_{\vert V} \subset 
Isom^+(V)$ of the action of the subgroup $H$ that they generate in $Diff^+(M)$ 
is finite. 

If the action on $V$ of the cyclic subgroups generated by two of the 
hyperelliptic rotations in $\Psi$ are conjugate in $H_{\vert V}$, one can 
conjugate the actions in $Diff^+(M)$ to coincide on $V$, since any 
diffeomorphism in $H_{\vert V}$ extends to $M$. Then by \cite[Lemma 10]{BoPa} 
these actions must coincide on $M$, contradicting the hypothesis that the 
conjugacy classes of cyclic subgroups generated by the hyperelliptic rotations 
in $\Psi$ are pairwise distinct in $Diff^+(M)$. Hence, the cyclic subgroup 
generated by the $k \geq 7$ hyperelliptic rotations in $\Psi$ are pairwise not 
conjugate in the finite group $H_{\vert V} \subset Isom^+(V)$.

Since the dual graph of the JSJ-decomposition of $M$ is a tree, a boundary 
torus $T \subset \partial V$ is separating and bounds a component $U_T$ of 
$M \setminus int(V)$. Since, by hypothesis, $\Gamma_f$ consists of a single 
vertex, no boundary component $T$ is setwise fixed by the finite group 
$H_{\vert V}$. This means that there is a hyperelliptic rotation 
$\psi_i \in \Psi$ of odd prime order $p_i$ such that the orbit of $U_T$ under 
$\psi_i$ is the disjoint union of $p_i$ copies of $U_T$. In particular $U_T$ 
projects homeomorphically onto a knot exterior in the quotient 
$\S^3 = M/\psi_i$. Therefore on each boundary torus $T = \partial U_{T} \subset 
\partial V$, there is a simple closed curve $\lambda_T$, unique up to isotopy, 
that bounds a properly embedded incompressible and $\partial$-incompressible 
Seifert surface $S_T$ in the knot exterior $U_T$. 

By pinching the surface $S_T$ onto a disc $\D^2$, in each component $U_T$ of 
$M \setminus int(V)$, one defines a degree-one map $\pi: M \longrightarrow M'$, 
where $M'$ is the rational homology sphere obtained by Dehn filling each 
boundary torus $T \subset V$ along the curve $\lambda_T$. 

For each hyperelliptic rotation $\psi_i$ in $\Psi$, of odd prime order $p_i$, 
the $\psi_i$-orbit of each component $U_T$ of $M \setminus int(V)$ consists of 
either one or $p_i$ elements. As a consequence, by \cite{Sa}, $\psi_i$ acts 
equivariantly on the set of isotopy classes of curves 
$\lambda_T \subset \partial V$. Hence, each $\psi_i$ extends 
to periodic diffeomorphism $\psi'_i$ of order $p_i$ on $M'$. Moreover, $M'$ is 
a $\ZZ/p_i$-homology sphere, since so is $M$ and $\pi: M \longrightarrow M'$ is 
a degree-one map. According to Smith theory, if $Fix(\psi')$ is non-empty on 
$M'$, then $\psi'_i$ is a rotation on $M'$. To see that 
$Fix(\psi')\neq\emptyset$ on $M'$ it suffices to observe that either
$Fix(\psi)\subset V$ or $\psi_i$ is a rotation of some $U_T$; in this latter
case, $\psi'_i$ must have a fixed point on the disc $\D^2$ onto which the
surface $S_T$ is pinched. To show that $\psi'_i$ is hyperelliptic it remains to 
show that the quotient $M'/\psi'_i$ is homeomorphic to $\S^3$.

Since $\psi_i$ acts equivariantly on the components $U_T$ of 
$M \setminus int(V)$ and on the set of isotopy classes of curves 
$\lambda_T \subset \partial V$, the quotient $\S^3 = M/\psi_i$ is obtained from 
the compact $3$-manifold $V/\psi_i$ by gluing knot exteriors (maybe solid tori) 
to its boundary components, in such a way that the boundaries of the Seifert 
surfaces of the knot exteriors are glued to the curves 
$\lambda_T/\psi_i \subset \partial V/\psi_i$.

In the same way, the rotation $\psi'_i$ acts equivariantly on the components 
$M '\setminus int(V)$ and on the set of isotopy classes of curves 
$\lambda_T \subset \partial V$. By construction, these components are solid 
tori, and either the axis of the rotation is contained in $V$ or there exists a 
unique torus $T \in \partial V$ such that the solid torus glued to $T$ to 
obtain $M'$ contains the axis. In the latter case, by \cite[Cor. 2.2]{EL}, the 
rotation $\psi'_i$ preserves a meridian disc of this solid torus and its axis 
is a core of $W_T$. It follows that the images in the quotient $M'/\psi'_i$ of 
the the solid tori glued to $\partial V$ are again solid tori. Hence 
$M'/\psi'_i$ is obtained from $\S^3$ by replacing each components of 
$\S^3 \setminus V/\psi_i$ by a solid torus, in such way that boundaries of 
meridian discs of the solid tori are glued to the curves 
$\lambda_T/\psi'_i \subset \partial V/\psi'_i$. It follows that $M'/\psi'_i$ is 
again $\S^3$.

So far we have constructed a closed orientable $3$-manifold $M'$ with a finite 
subgroup of orientation preserving diffeomorphisms $H_V$ that contains at least 
seven conjugacy classes of hyperelliptic  subgroups  of odd prime orders. Theorem~\ref{th:main4} implies that $M'$ must be 
$\S^3$, and thus by the orbifold theorem \cite{BLP} $H_V$ is conjugate to 
a finite subgroup of $SO(4)$. In particular the subgroup $H_0 \subset H_V$ 
generated by the subset $\Psi_0$ of at least $4$ hyperelliptic rotations with 
pairwise distinct odd prime orders must be solvable. Therefore, by 
Proposition~\ref{pro:solvable} the induced rotations commute on $M'$ and, 
by restriction, the hyperelliptic rotations in $\Psi_0$ commute on $V$. 
\qed

\medskip

In the final step of the proof we extend the commutativity on $M_f$ to the 
whole manifold $M$. The proof of this step is analogous to the one given in 
\cite{BPZ}, since the proof there was not using the homology assumption. We 
give the argument for the sake of completeness.

\begin{Proposition}\label{prop:whole commuting}
The $k_0 \geq 4$ hyperelliptic rotations in $\Psi_0$ commute, up to conjugacy 
in $Diff^+(M)$, on $M$.
\end{Proposition}

\demo

Let $\Gamma_c$ be the largest subtree of $\Gamma$ containing $\Gamma_f$, such 
that, up to conjugacy in $Diff^{+}(M)$, the rotations in $\Psi_0$
commute on the corresponding invariant submanifold $M_c$ of $M$. We shall 
show that $\Gamma_c=\Gamma$. If this is not the case, we can choose an edge 
contained in $\Gamma$ corresponding to a boundary torus $T$ of $M_c$. Denote by 
$U_T$ the submanifold of $M$ adjacent to $T$ but not containing  $M_c$ and by 
$V_T \subset U_T$ the geometric piece adjacent to $T$.

Let $H_0 \subset Diff^+(M)$ be the group of diffeomorphisms of $M$ generated by 
the set of geometric hyperelliptic rotations 
$\Psi_0 =\{\psi_i, \, i = 1,\dots,k_0 \}$. Since the rotations in $\Psi_0$ 
commute on $M_c$ and have pairwise distinct odd prime orders, the restriction 
of $H_0$ on on $M_c$ is a cyclic group with order the product of the orders of 
the rotations. Since $\Gamma_f \subset\Gamma_c$, the $H_0$-orbit of $T$ cannot 
be reduced to only one element. Moreover each rotation $\psi \in \Psi_0$ either 
fixes $T$ or acts freely on the orbit of $T$ since its order is prime.

If no rotation in $\Psi_0$ leaves $T$ invariant, the $H_0$-orbit of $T$ 
contains as many elements as the product of the orders of the rotations, for 
they commute on $M_c$. In particular, only the identity (which extends to $U$) 
stabilises a torus in the $H_0$-orbit of $T$. Note that all components of 
$\partial M_c$ that are in the $H_0$-orbit of $T$ bound a manifold homeomorphic 
to $U_T$.

Since the rotation $\psi_i$ acts freely on the $H_0$-orbit of $U_T$, $U_T$ is a 
knot exterior in the quotient $M/\psi_i=\S^3$. Hence there is a well defined 
meridian-longitude system on $T = \partial U_T$ and also on each torus of the 
$H_0$-orbit of $T$. This set of meridian-longitude systems is cyclically 
permuted by each $\psi_i$ and thus equivariant under the action of $H_0$.

Let $M_c/H_0$ be the quotient of $M_c$ by the induced cyclic action of $H_0$ on 
$M_c$. Then there is a unique boundary component 
$T' \subset \partial (M_c/H_0)$ which is the image of the $H_0$-orbit of $T$. 
We can glue a copy of the knot exterior $U_T$ to $M_c/H_0$ along $T'$ by 
identifying the image of the meridian-longitude system on $\partial U_T$ with 
the projection on $T'$ of the equivariant meridian-longitude systems on the 
$H_0$-orbit of $T$. Denote by $N$ the resulting manifold. For all 
$i=1,\dots, k_0$, consider the cyclic (possibly branched) cover of $N$ of order 
$q_i=\prod_{j\neq i} p_j$ which is induced by the cover 
$\pi_i:M_c/\psi_i\longrightarrow M_c/H_0$. This makes sense because  
$\pi_1(T')\subset{\pi_i}_*(\pi_1(M_c/\psi_i))$. Call $\tilde N_i$ the total 
space of such covering. By construction it follows that $\tilde N_i$ is the 
quotient $(M_c\cup H_{0}\cdot U_T)/\psi_i$. This implies that the $\psi_i$'s  
commute on $M_c\cup H_{0}\cdot U_T$ contradicting the maximality of $\Gamma_c$. 

We can thus assume that some rotations fix $T$ and some do not. Since all 
rotations commute on $M_c$ and have pairwise distinct odd prime orders, we see 
that the orbit of $T$ consists of as many elements as the product of the orders 
of the rotations which do not fix $T$ and each element of the orbit is fixed by 
the rotations which leave $T$ invariant. The rotations which fix $T$ commute on 
the orbit of $V_T$ according to Claim~\ref{claim:invariant edge} and Lemma 
\ref{lem:gluing}, and form a cyclic group generated by, say, $\gamma$. The 
argument for the previous case shows that the rotations acting freely on the 
orbit of $T$ commute on the orbit of $U_T$ and thus on the orbit of $V_T$, and form 
again a cyclic group generated by, say, $\eta$. To reach a contradiction to the 
maximality of $M_c$, we shall show that $\gamma$, after perhaps some 
conjugacy, commutes with $\eta$ on the $H_0$-orbit of $V_T$, in other words 
that $\gamma$ and $\eta\gamma\eta^{-1}$ coincide on $H_{0}\cdot V_T$. Since 
$\eta$ acts freely and transitively on the $H_0$-orbit of $V_T$ there is a 
natural and well-defined way to identify each element of the orbit 
$H_{0}\cdot V_T$ to $V_T$ itself. Note that this is easily seen to be the case if
$V_T$ is hyperbolic: this follows from Claim~\ref{claim:invariant edge} and
Claim~\ref{claim:slope}. We now consider the case when $V_T$ is Seifert fibred.

\begin{Claim}\label{claim:consistent action}
Assume that $V_T$ is Seifert fibred and that the restriction of $\gamma$ 
induces a non-trivial action on the base of $V_T$. Then $\gamma$ induces a 
non-trivial action on the base of each component of the $H_0$-orbit of $V_T$. 
Moreover, up to conjugacy on $H_{0}\cdot V_T \setminus V_T$ by diffeomorphisms 
which are the identity on $H_{0}\cdot T$ and extend to $M$ , we can assume that 
the restrictions of $\gamma$ to these components induce the same permutation of 
their boundary components and the same action on their bases.
\end{Claim}

\demo

By hypothesis $\gamma$ and $\eta\gamma\eta^{-1}$ coincide on $\partial M_c$. 
The action of $\gamma$ on the base of $V_T$ is non-trivial if and only if its 
restriction to the boundary circle of the base corresponding to the fibres of 
the torus $T$ is non-trivial. Therefore the action of $\gamma$ is non-trivial 
on the base of each component of $H_{0}\cdot V_T$.
 
By Lemma~\ref{lem:seifert} and taking into account the fact that $V_T$ is a
geometric piece in the JSJ-decomposition of the knot exterior $U_T$, the only 
situation in which the action of $\gamma$ on the base of $V_T$ is non-trivial 
is when the base of $V_T$ consists of a disc with $p$ holes, where $p$ is the 
order of one of the rotations that generate $\gamma$.
%, the only one whose action is non-trivial on the base of the fibration. 
Moreover, the restriction 
of $\gamma$ to the elements of $H_{0}\cdot V_T$ cyclically permutes their 
boundary components which are not adjacent to $M_c$. Up to performing Dehn 
twists, along vertical tori inside the components of $H_{0}\cdot V_T \setminus 
V_T$, which permute the boundary components, we can assume that the restriction 
of $\gamma$ induces the same cyclic permutations on the boundary components of 
each element of the orbit $H_{0}\cdot V_T$. We only need to check that Dehn 
twists permuting two boundary components extend to the whole manifold $M$. This 
follows from the fact that the manifolds adjacent to these components are all 
homeomorphic and that Dehn twist act trivially on the homology of the boundary. 

Since the actions of the restrictions of $\gamma$ on the bases of the elements 
of $H_{0}\cdot V_T$ are combinatorially equivalent, after perhaps a further 
conjugacy by an isotopy, the different restrictions can be chosen to coincide 
on the bases.
\qed

\medskip

By Claim~\ref{claim:invariant edge} and Claim~\ref{claim:consistent action} we 
can now deduce that the restrictions of $\gamma$ and $\eta\gamma\eta^{-1}$ 
to the $H_0$-orbit of $V_T$ commute, up to conjugacy of $\gamma$ which is the 
identity on the $H_0$-orbit of $T$. Since $\gamma$ and $\eta\gamma\eta^{-1}$ 
coincide on this $H_0$-orbit of $T$, we can conclude that they coincide on the 
$H_0$-orbit of $V_T$. This finishes the proof of 
Proposition~\ref{prop:whole commuting}.
\qed
\smallskip

Since there are at least four hyperelliptic rotations with paiwise distinct odd 
prime orders in $\Psi_0$, Proposition~\ref{prop:finitereduction} is consequence 
of Proposition \ref{prop:whole commuting} and 
Proposition~\ref{prop:commuting rotations}, like in the solvable case. 
\qed

\begin{remark}
As we have seen, the strategy to prove that an irreducible manifold $M$ with 
non-trivial JSJ-decomposition cannot admit more than six conjugacy classes of
hyperelliptic subgroups of odd prime order inside
$Diff^{+}(M)$ consists in modifying by conjugacy any given set of hyperelliptic
rotations so that the new hyperelliptic rotations commute pairwise. Note that
this strategy cannot be carried out in general if the orders are not pairwise
coprime (see, for instance, \cite[Section 4.1]{BoPa} where the case of two 
hyperelliptic rotations of the same odd prime order, generating non conjugate 
subgroups, is considered). Similarly, for hyperelliptic rotations of arbitrary
orders $>2$ lack of commutativity might arise on the Seifert fibred pieces of
the decomposition, as it does for the Seifert fibred manifolds constructed in
Proposition~\ref{prop:circle bundles}. 
\end{remark}

%%%%%%%%%%%%%%%%%%%%%%%%%%%%%%%%%%%%%%%%%%%%%%%%%%%%%%%%%%%%%%%%%%%%%%%%%%%%

\section{Appendix: non-free finite group actions on rational homolgy spheres}
\label{sec:nonfree}

In this section we will show that every finite group $G$ admits a faithful 
action by orientation preserving diffeomorphisms on some rational homology 
sphere so that some elements of $G$ have non-empty fixed-point sets.

Cooper and Long's construction of $G$-actions on rational homolgy spheres in 
\cite{CL} consists in starting with a $G$-action on some $3$-manifold and then
modifying the original manifold, notably by Dehn surgery, so that the new 
manifold inherits a $G$-action but has ``smaller" rational homology. Since 
their construction does not require that the $G$-action is free, it can be used 
to prove the existence of non-free $G$-actions. We will thus start by
exhibiting non-free $G$-actions on some $3$-manifold before pointing out what
need to be taken into account in this setting in order for Cooper and Long's 
construction to work.

Since every (non-trivial) cyclic group acts as a group of rotations of $\S^3$, for simplicity
we will assume that $G$ is a finite non-cyclic group.

\begin{Claim}\label{claim:non-free action}
Let $G$ be a finite non-cyclic group. There is a closed, connected, orientable 
$3$-manifold $M$ on which $G$ acts faitfully, by orientation preserving
diffeomorphisms so that there are $g\in G\setminus\{1\}$ with the property that
$Fix(g)\neq\emptyset$.
\end{Claim}

\demo

Let $k\ge 2$ and let  $\{g_i\}_{1\le i\le k+1}$ be a system of generators for 
$G$ satisfying the following requirements:
\begin{itemize}
\item for all $1\le i\le k+1$, the order of $g_i$ is $n_i\ge 2$;
\item $g_{k+1}=g_1g_2\dots g_k$.
\end{itemize}
Since $G$ is not cyclic these conditions are not difficult to fulfill, and
actually they can be fulfilled even in the case when $G$ is cyclic for an
appropriate choice of the set $\{g_i\}_{1\le i\le k+1}$.

Consider the free group of rank $k$ that we wish to see as the fundamental
group of a $(k+1)$-punctured $2$-sphere: each generator $x_i$ corresponds to a 
loop around a puncture of the sphere so that a loop around the $k+1$st puncture 
corresponds to the element $x_1x_2\dots x_k$. Build an orbifold $\O$ by 
compactifying the punctured-sphere with cone points so that the $i$th puncture
becomes a cone point of order $n_i$. The resulting orbifold has (orbifold)
fundamental group with the following presentation:
$$\langle x_1, x_2, \dots, x_k, x_{k+1}\mid \{ x_i^{n_i}\}_{1\le i\le k+1},\, 
x_1\ldots x_kx_{k+1}^{-1} \rangle.$$
Clearly this group surjects onto $G$. Such surjection gives rise to an orbifold
covering $\Sigma\longrightarrow\O$, where $\Sigma$ is an orientable surface on
which $G$ acts in such a way that each element $g_i$ has non-empty fixed-point
set. One can consider the $3$-manifold $\Sigma \times \S^1$: the action of $G$
on $\Sigma$ extends to a product action of $G$ on $\Sigma \times \S^1$ which is
trivial on the $\S^1$ factor.
\qed 
\smallskip

To be able to repeat Cooper and Long's construction it is now easy to observe
that it is always possible to choose $G$-equivariant families of simple closed
curves in $M$ so that they miss the fixed-point sets of elements of $G$ and
either their homology classes generate $H_1(M;\QQ)$ (so that the hypothesis of 
\cite[Lemma 2.3]{CL} are fulfilled when we choose $X$ to be the exterior of
such families) or the family is the $G$-orbit of a representative of some
prescribed homology class (as in the proof of \cite[Proposition 2.5]{CL}).

\paragraph{Acknowledgement} 
The authors are indebted to M. Brou\'e for valuable discussions on the topics
of the paper.

\textsc{Aix-Marseille Universit\'e, CNRS, Centrale Marseille, I2M, UMR 7373,
13453 Marseille, France}

{michel.boileau@univ-amu.fr}

\textsc{Dipartimento di Matematica e Fisica ``Niccol\`o Tartaglia'', 
Universit\`a Cattolica del Sacro Cuore, Via Musei 41, 25121 Brescia, Italy}

{clara.franchi@unicatt.it}

\textsc{Dipartimento di Matematica e Geoscienze, Universit\'a degli studi di 
Trieste, Via Valerio 12/1, 34127 Trieste}

{mmecchia@units.it}

\textsc{Aix-Marseille Universit\'e, CNRS, Centrale Marseille, I2M, UMR 7373,
13453 Marseille, France}

{luisa.paoluzzi@univ-amu.fr}

\textsc{Dipartimento di Matematica e Geoscienze, Universit\'a degli studi di 
Trieste, Via Valerio 12/1, 34127 Trieste}

{zimmer@units.it}


\footnotesize

\begin{thebibliography}{XXXXX}



\bibitem[A]{A} M. Aschbacher,
\emph{Finite groups theory}, 
Cambridge studies in advanced mathematics \textbf{10}, 
Cambridge University Press, 1986.

%\bibitem[B3MP]{B3MP} L. Bessi\`eres, G. Besson, M. Boileau, S. Maillot, and 
%J. Porti,
%\emph{Geometrisation of $3$-manifolds},
%Tracts in Math. \textbf{13}
%European Mathematical Society, 2010.

\bibitem[BMP]{BMP} M. Boileau, S. Maillot, and J. Porti, 
\emph{Three-dimensional orbifolds and their geometric structures}, 
Panoramas et Synth\`eses \textbf{15}, 
Soci\'et\'e Math\'ematique de France, Paris, 2003.


\bibitem[BLP]{BLP} M. Boileau, B. Leeb and J. Porti, 
\emph{ Geometrization of $3$-dimensional orbifolds},
Annals of Math. \textbf{162},  (2005), 195-290.

\bibitem[BOt]{BOt}  M. Boileau and J.P. Otal.
\emph{ Scindements de Heegaard et groupe des hom\'eotopies des petites vari\'et\'es de Seifert},
Invent. Math. \textbf{106}, (1991), 85-107.

\bibitem[BoPa]{BoPa} M. Boileau and L. Paoluzzi, 
\emph{On cyclic branched coverings of prime knots}, 
J. Topol. \textbf{1}, (2008), 557-583.

\bibitem[BPZ]{BPZ} M. Boileau, L. Paoluzzi and B. Zimmermann, 
\emph{A characterisation of ${\bf S}^3$ among homology spheres}, 
The Zieschang Gedenkschrift. Geom. Topol. Monogr. \textbf{14}, 
Geom. Topol. Publ., Coventry, (2008), 83-103.

\bibitem[BoP]{BoP} M. Boileau and J. Porti, 
\emph{Geometrization of 3-orbifolds of cyclic type}, 
Ast\'erisque Monograph, \textbf{272}, 2001.


\bibitem[BS]{BS} F. Bonahon and L. Siebenmann, 
\emph{The characteristic splitting of irreducible compact 3-orbifolds}, 
Math. Math. \textbf{278}, (1987), 441-479.


\bibitem[BZH]{BZH}
G. Burde and H. Zieschang, M. Heusener,
{\em Knots}, Third edition, De Gruyter Studies in Mathematics, 5. De Gruyter, 
Berlin, 2014.



%NOT CITED ANYMORE/ANYWHERE: REMOVE???
%\bibitem[BM]{BM} M. Brou\'e and G. Malle, 
%\emph{Th\'eor\`emes de Sylow g\'en\'eriques pour les groupes r\'eductifs sur 
%les corps finis}, Math. Ann. \textbf{292}, (1992), 241-262.


\bibitem[CL]{CL} D. Cooper and D. D. Long, 
\emph{Free actions of finite groups on rational homology $3$-spheres}, 
Topology Appl. \textbf{101}, (2000), 143–148.

%NOT CITED ANYMORE/ANYWHERE: REMOVE???
%\bibitem[Co]{Co} J. H. Conway, R. T. Curtis, S. P. Norton, R. A. Parker, and  
%R. A. Wilson, 
%\emph{Atlas of Finite Groups}, 
%Oxford University Press, 1985.


\bibitem[CHK]{CHK} D. Cooper, C. Hodgson, and S. Kerckhoff,
\emph{Three-dimensional orbifolds and cone-manifolds}, 
MSJ Memoirs \textbf{5}, 2000.

\bibitem[DL]{DL} J. Dinkelbach and B. Leeb, 
\emph{Equivariant Ricci flow with symmetry and applications to finite group actions on 3-manifolds}, Geom. Topol. \textbf{13}, (2009), 1129–1173. 


\bibitem[Dun]{Dun} W. D. Dunbar, 
\emph{Geometric Orbifolds}, 
Rev. Mat. Univ. Comp. Madrid, \textbf{1}, 1998, 67-99.

\bibitem[EL]{EL} A. L. Edmonds and C. Livingston, 
\emph{Group actions on fibered three-manifolds}, 
Comm. Math. Helv. \textbf{58}, (1983), 529-542.

\bibitem[Gr]{Gr} J. E. Greene, 
\emph{Lattices, graphs, and Conway mutation.}
Invent. Math. \textbf{192} (2013), no. 3, 717–750. 


\bibitem[Go]{Go} C. McA. Gordon, 
\emph{Some aspects of classical knot theory}, 
Knot Theory, Proceedings, Plans-sur-Bex, Switzerland (J.C. Hausmann ed.), 
Lect. Notes Math. \textbf{685}, (1977), Springer-Verlag, 1-60.

%\bibitem[G]{G} D. Gorenstein, 
%\emph{The classification of finite simple groups. Volume 1: Groups of
%Noncharacteristic 2 Type} 
%The University Series in Mathematics, \textbf{1}, 
%Plenum Press, New York, 1983.

%NOT CITED ANYMORE/ANYWHERE: REMOVE???
%\bibitem[GLS1]{GLS1} D. Gorenstein, R. Lyons, and R. Solomon,
%\emph{The classification of the finite simple groups, Number 1},
%Math. Surveys Monogr. \textbf{40.1}, 
%Amer. Math. Soc., Providence, RI, 1995.

%NOT CITED ANYMORE/ANYWHERE: REMOVE???
%\bibitem[GLS2]{GLS2} D. Gorenstein, R. Lyons, and R. Solomon,
%\emph{The classification of the finite simple groups, Number 2},
%Math. Surveys Monogr. \textbf{40.2}, 
%Amer. Math. Soc., Providence, RI, 1996.

\bibitem[GLS]{GLS} D. Gorenstein, R. Lyons, and R. Solomon,
\emph{The classification of the finite simple groups, Number 3},
Math. Surveys Monogr. \textbf{40.3}, 
Amer. Math. Soc., Providence, RI, 1998.

\bibitem[Hil]{Hil} J. Hillman, 
\emph{Links with infinitely many semifree periods are trivial}, 
Arch. Math. \textbf{42}, (1984), 568-572.


\bibitem[H]{H} B. Huppert, 
\emph{Endlichen Gruppen I}, 
Springer-Verlag, New York, 1968. 

%\bibitem[K]{K}  E. I. Khukhro, 
%\emph{Nilpotent Groups and Their Automorphisms}, 
%de Gruyter-Verlag, Berlin, 1993. 


\bibitem[KL]{KL} B. Kleiner and J. Lott.
\emph{Local Collapsing, Orbifolds, and Geometrization}, 
Ast\'erisque Monograph \textbf{365}, 2014.


\bibitem[JS]{JS} W. H. Jaco and P. B. Shalen, 
\emph{Seifert fibred spaces in 3-manifolds}, 
Mem. Amer. Math. Soc. \textbf{220}, 1979.

\bibitem[J]{J} K. Johannson, 
\emph{Homotopy equivalence of 3-manifolds with boundary}, 
Lecture Notes in Math. \textbf{761}, Springer-Verlag, Berlin, 1979.

\bibitem[Ka]{Ka} A. Kawauchi,
\emph{Topological imitations and Reni-Mecchia-Zimmermann's conjecture},
Kyungpook Math. J. \textbf{46} (2006), 1–9. 

\bibitem[Ko]{Ko} S. Kojima, 
\emph{Determining knots by branched covers}, 
in Low Dimensional Topology and Kleinian groups, 
London Math. Soc. Lecture Note Ser. \textbf{112}, 
Cambridge Univ. Press (1986), 193-207.



\bibitem[Mec1]{Mec0} M. Mecchia,
\emph{How hyperbolic knots with homeomorphic cyclic branched coverings are 
related}, 
Topology Appl. \textbf{121}, (2002), 521-533.


\bibitem[Mec2]{Mec} M. Mecchia,
\emph{Finite groups acting on 3-manifolds and cyclic branched coverings of
knots}, 
The Zieschang Gedenkschrift. Geom. Topol. Monogr. \textbf{14}, 
Geom. Topol. Publ., Coventry, (2008), 393-416.

\bibitem[MR]{MR} M. Mecchia and M. Reni,
\emph{Hyperbolic 2-fold branched coverings of links and their quotients},
Pacific J. Math. \textbf{202}, (2002), 429–447. 

%\bibitem[MZ]{MZ} M. Mecchia and B. Zimmermann, 
%\emph{The number of knots and links with the same 2-fold branched covering}, 
%Quart. J. Math. \textbf{55} (2004) 69–76.


\bibitem[MZ]{MZ} M. Mecchia and B. Zimmermann, 
\emph{On finite groups acting on $\ZZ_2$-homology $3$-spheres}, 
 Math. Z. \textbf{248}, (2004), 675-693.


\bibitem[MSY]{MSY} W. H. Meeks, L. Simon, and S.-T. Yau, 
\emph{Embedded minimal surfaces, exotic spheres, and manifolds with positive Ricci curvature},
Ann. of Math. \textbf{116} (1983), 621–659.


\bibitem[MS]{MS} W. H. Meeks and P. Scott,
\emph{Finite group actions on 3-manifolds},
Invent. Math. \textbf{86}, (1986), 287-346.

\bibitem[Mon1]{Mon1} J. M. Montesinos, 
\emph{Variedades de Seifert que son recubridores ciclicos ramificados de dos 
hojas},  
Bol. Soc. Mat. Mexicana \textbf{18}, (1973), 1-32.

\bibitem[Mon2]{Mon2} J. M. Montesinos, 
\emph{Rev\^etements ramifi\'es de noeuds, espaces fibr\'es de Seifert et 
scindements de Heegaard},
Publicaciones del Seminario Mathematico Garcia de Galdeano, Serie II, Seccion 
3, (1984).

\bibitem[MonW]{MonW} J. M. Montesinos and W. Whitten
\emph{Constructions of two-fold branched covering spaces.} Pac. J. Math. \textbf{125}, 415–446 (1986) 


\bibitem[NS]{NS} W. Neumann and G. Swarup, 
\emph{Canonical decompositions of $3$-manifolds}, 
Geom. Topol. \textbf{1}, (1998), 21-40. 

\bibitem[Or]{Or} P. Orlik, 
\emph{Seifert manifolds}, 
Lecture Notes in Mathematics \textbf{291}, Springer, 1972.

\bibitem[P]{P} L. Paoluzzi,
\emph{On hyperbolic type involutions},
Rend. Istit. Mat. Univ. Trieste \textbf{32} suppl. 1, (2001), 257-288.

\bibitem[Re]{Re} M. Reni, 
\emph{On $\pi$-hyperbolic knots with the same 2-fold branched coverings}, 
Math. Ann. \textbf{316}, (2000), no. 4, 681–697

%\bibitem[RZ]{RZ} M. Reni and B. Zimmermann, 
%\emph{Finite simple groups acting on 3-manifolds and homology spheres}, 
%Rend. Istit. Mat. Univ. Trieste \textbf{32} suppl. 1, (2001), 300-313.


\bibitem[RZ]{RZ} M. Reni and B. Zimmermann, 
\emph{Hyperbolic 3-manifolds as cyclic branched coverings}, 
Comment. Math. Helv. \textbf{76}, (2001), 300-313.

\bibitem[Sa]{Sa} M. Sakuma, 
\emph{Uniqueness of symmetries of knots},
Math. Z. \textbf{192}, (1986), 225-242.

%\bibitem[Sco]{Sco} P. Scott, 
%\emph{The geometries of 3-manifolds}, 
%Bull. London Math. Soc. \textbf{15}, (1983), 401-487.

\bibitem[Sco]{Sco} P. Scott,
\emph{Homotopy implies isotopy for some Seifert fibre spaces},
Topology \textbf{24}, (1985), 341-351.



%\bibitem[Su1]{Su1} M. Suzuki, 
%\emph{Group theory I}, 
%Grundlehren Math. Wiss. \textbf{247}, Springer-Verlag, Berlin, 1982.

\bibitem[Su]{Su} M. Suzuki, 
\emph{Group theory II}, 
Grundlehren Math. Wiss. \textbf{248}, Springer-Verlag, New York, 1982.


\bibitem[V]{V}  O. Ja. Viro 
\emph{Nonprojecting isotopies and knots with homeomorphic coverings}, 
Zap. Nauč. Semin. POMI  \textbf{66}, 133–147, 207–208, Studies in topology, II (1976) 

\bibitem[Wa]{Wa} F. Waldhausen, 
\emph{Eine Klasse von $3$-dimensionalen Mannigfaltigkeiten. I, II}, 
Invent. Math. \textbf{3}, (1967), 308–333; ibid. \textbf{4} (1967), 87–117.




\bibitem[Z1]{Z1} B. Zimmermann, 
\emph{On the Hantzsche-Wendt Manifold}, 
Mh. Math. \textbf{110}, (1990) , 321-327.



\bibitem[Z2]{Z2} B. Zimmermann, 
\emph{On hyperbolic knots with homeomorphic cyclic branched coverings}, 
Math. Ann. \textbf{311}, (1998), 665-673.


%\bibitem[Ko1]{Ko1} S. Kojima, 
%\emph{Bounding finite group acting on 3-manifolds}, 
%Math. Proc. Camb. Phil. Soc. \textbf{96}, (1984), 269-281.
%\bibitem[Ko3]{Ko3} S. Kojima, 
%\emph{Isometry transformations of hyperbolic 3-manifolds}
￼%Topology and its Applications \textbf{29} (1988) 297-307
%\bibitem[LMM]{LMM} A. Lucchini, F. Menegazzo and M Morigi,
%\emph{On the existence of a complement for a finite simple group in its
%automorphism group}, Illinois J. Math. \textbf{47}, (2003), 395-418.
%NOT CITED ANYMORE/ANYWHERE: REMOVE???
%\bibitem[MaT]{MaT} G. Malle and D. Testerman 
%\emph{Linear algebraic groups and finite groups of Lie type.} 
%Cambridge Studies in Advanced Mathematics \textbf{133}, 
%Cambridge University Press, Cambridge, 2011.
%\bibitem[Macc]{Macc} A. M. MacBeath, \emph{Generators of the linear fractional
%groups}, Proc. Symp. Pure Math. \textbf{12}, (1968), 14-32.


\end{thebibliography}
\end{document}